\documentclass[article]{siamart}

\usepackage{amsmath,amssymb}
\usepackage{array}
\usepackage{caption}
\usepackage{hyperref}
\usepackage[utf8]{inputenc}
\usepackage{mathtools}
\usepackage{color,graphicx}
\usepackage{thmtools, thm-restate}
\usepackage{enumitem}
\usepackage{algorithmicx}
\usepackage{algpseudocode}
\usepackage{tikz}
\usetikzlibrary{arrows}
\usetikzlibrary{3d,calc}

\newcommand{\cH}{\mathcal{H}}

\newcommand{\cO}{\mathcal{O}}

\newcommand{\cP}{\mathcal{P}}
\newcommand{\cR}{\mathcal{R}}
\newcommand{\cS}{\mathcal{S}}

\newcommand{\beq}{\begin{equation}}
\newcommand{\eeq}{\end{equation}}
\def\bals#1\eals{\begin{align*} #1 \end{align*}}
\def\bal#1\eal{\begin{align} #1 \end{align}}

\newcommand\Dom\Omega

\newcommand\RR{\mathbb{R}}
\newcommand\PP{\mathbb{P}}

\newcommand\ZZ{\mathbb{Z}}

\newcommand\Lap\Delta

\newcommand\abs[1]{\left\lvert #1 \right\rvert}

\def\bpde#1\epde{\[\left\{\begin{aligned}#1\end{aligned}\right. \]}
\def\inbpde#1\inepde{\left\{\begin{aligned}#1\end{aligned}\right.}
\def\binpde#1\einpde{\left\{\begin{aligned}#1\end{aligned}\right.}

\newcommand\prn[1]{\left( {#1} \right)}

\newcommand\Norm[2]{\left\lVert { #1 } \right\rVert_{#2}}

\def\grad{\nabla}

\def\cB{\mathcal{B}}

\def\cK{\mathcal{K}}

\def\cR{\mathcal{R}}
\def\half{\frac{1}{2}}

\def\p{\partial}

\def\b0{\mathbf{0}}

\def\bbmat{\begin{bmatrix}[r]}
\def\ebmat{\end{bmatrix}}

\newcommand{\barr}{\begin{array}}
\newcommand{\ea}{\end{array}}
\newcommand{\bea}{\begin{eqnarray}}
\newcommand{\eea}{\end{eqnarray}}
\newcommand{\bt}{\begin{table}}
\newcommand{\et}{\end{table}}

\DeclareMathOperator\argmin{argmin}

\theoremstyle{plain}
\theoremstyle{definition}

\numberwithin{equation}{section}

\begin{document}
  
\title{Dimensional splitting of hyperbolic partial differential equations 
       using the Radon transform}

\author{Donsub Rim%
  \thanks{Department of Applied Mathematics, University of Washington, Seattle
WA 98195 (\email{drim@uw.edu}). %
   Current affiliation: Department of Applied Physics and Applied Mathematics, %
   Columbia University, New York, NY 10027 (\email{dr2965@columbia.edu}).}
}

\maketitle

\begin{abstract} 
We introduce a dimensional splitting method based on the intertwining property
of the Radon transform, with a particular focus on its applications related to
hyperbolic partial differential equations (PDEs).  This dimensional splitting
has remarkable properties that makes it useful in a variety of contexts,
including multi-dimensional extension of large time-step (LTS) methods,
absorbing boundary conditions, displacement interpolation, and
multi-dimensional generalization of transport reversal \cite{reversal}. 
\end{abstract}

\section{Introduction}

Dimensional splitting provides the simplest approach to obtaining a
multi-dimensional method from a one-dimensional method
\cite{acta_splitting,strangsplit,Crandall1980,fvmbook}.  Although extremely
powerful, existing splitting methods do not preserve a special feature that is
easily obtained for 1D methods. For 1D hyperbolic partial differential
equations (PDEs) of the type
\begin{equation}
q_t + A q_x  = 0
\end{equation}
where $A$ is a constant diagonalizable matrix with real and distinct
eigenvalues, one can devise large time-step (LTS) methods that allow the
solution to be solved up to any time without incurring excessive numerical
diffusion \cite{largetimestep1, largetimestep2, largetimestep3}.  Previous
splitting methods do not lead to such LTS methods in multi-dimensions. 

In this paper, we introduce a dimensional splitting method that allows
multi-dimensional linear constant coefficient hyperbolic problems to be solved
up to desired time. The method relies on the intertwining property of Radon
transforms \cite{Helgason2011,ctbook}, thereby transforming a multi-dimensional
problem into a family of one-dimensional ones.  Simply by applying an 1D LTS
method on each of these one-dimensional problems, one obtains a
multi-dimensional LTS method.  While this intertwining property is well-known
and is utilized to analyze PDEs in standard references \cite{laxphilips},
it has not been used for constructing multi-dimensional numerical methods, to
the best of our knowledge.

The method also has implications for the problem of imposing absorbing boundary
conditions, a problem that has received sustained interest over many decades
\cite{abc,pml,stretched-coords,pml3d}.  By using the Radon transform, the
splitting decomposes multi-dimensional waves into planar ones, thereby allowing
a separate treatment of each incident planar wave near the boundary. This
yields the desired absorbing boundary conditions in odd dimensions, and in even
dimensions one obtains an approximation up to $\cO(1/t)$ that does not cause
spurious reflections.

Another useful application of this dimensional splitting is in displacement
interpolation, a concept that arises naturally in optimal transport
\cite{villani2008optimal}. Our interest in displacement interpolation is
motivated by model reduction. To construct reduced order models for typical
hyperbolic problems, one cannot rely solely on linear subspaces
\cite{amsallem,Carlberg15}, and it is necessary to interpolate over the
Lagragian action \cite{marsden1, marsden2, Schulze15, reversal}. In a single
spatial dimension this can be done in a relatively straightforward manner,
owing to the LTS methods available for 1D \cite{reversal}.  The
multi-dimensional LTS method is useful also for the multi-dimensional extension
of displacement interpolation, and this in turn will yield a straightforward
way for low-dimensional information to be extracted for multi-dimensional
hyperbolic problems.

For the dimensional splitting to be computationally successful, one requires an
algorithm for computing the Radon transform and its inverse efficiently.
Throughout this paper we use the approximate discrete Radon transform (ADRT),
also called simply the discrete Radon transform (DRT), devised in
\cite{brady,GD96}. We will refer to ADRT as DRT.  It is a fast algorithm with
the computational cost of  $\cO(N^2 \log N)$ for an $N \times N$ image or
grid\footnote{The term \emph{grid} (\emph{cell}) is a more appropriate term for
our PDE applications, but DRT originally comes from imaging literature so we
will sometimes also use the term \emph{image} (\emph{pixel}),
interchangeably.}, and the efficiency is obtained through a geometric
recursion of so-called \emph{digital lines}.  The inversion algorithm using the
full multi-grid method appeared in \cite{press}, but here we adopt a simpler
approach by making use of the conjugate gradient algorithm \cite{greenbaum} for
the inversion. 

This paper is organized as follows. In Section \ref{sec:radonsplit}, we give a
review of the intertwining property of the Radon transform and introduce the
dimensional splitting method. In Section \ref{sec:drt}, we give a brief
introduction to the DRT algorithm and discuss its inversion. In Section
\ref{sec:radonsplituse}, we discuss its applications in absorbing boundary and
in displacement interpolation.  In this paper we will fully implement only
constant coefficient linear problems in spatial dimension two, although we will
also discuss how the splitting can be extended to fully nonlinear problems and
to higher spatial dimensions.  Further investigations into these and other
related topics will be mentioned in Section \ref{sec:conclusion}.

The Radon transform was introduced by Johann Radon \cite{radon} and has been a
major subject of study, primarily due to its use in medical imaging but also as
a general mathematical and computational tool.

\section{Dimensional splitting using the Radon transform}\label{sec:radonsplit}

In this section, we briefly review the intertwining property of the Radon
transform \cite{Helgason2011,ctbook} then show that it can be used as a
dimensional splitting tool that extends the large time-step (LTS) operator to
multiple spatial dimensions. It preserves the ability to take large time-steps
without loss of accuracy in the constant coefficient case. Moreover, this
splitting can potentially be used for fully nonlinear problems as well, in
a similar manner to the other splitting methods, with the usual CFL condition
for the time-step.

\subsection{Intertwining property of Radon transforms}

The \emph{Radon transform}  $\hat{\varphi}: S^{n-1} \times \RR \to \RR$ of the
function $\varphi: \RR^n \to \RR$ is defined as 
\begin{equation}	
 \hat{\varphi}(\omega, s) = \cR \varphi (\omega, s)
= \int_{x \cdot \omega = s} 
        \varphi(x) \, \textrm{d}m(x),
  \label{eq:radon}
\end{equation}
in which $\mathrm{d}m$ is the Euclidean measure over the hyperplane.  For any
fixed pair $(\omega, s) \in S^{n-1} \times \RR$, the set $\{x = (x_1,
x_2, ..., x_n) \in \RR^n : x \cdot \omega  = s\}$ defines a hyperplane, so
the transform is simply an integration of the function over this hyperplane.
We will denote the space of hyperplanes parametrized above by $\PP^n$.
In effect, $\hat{\varphi}$ decomposes $\varphi$ into planar waves in the
direction of $\omega$.

The \emph{back-projection} is defined as the dual of $\cR$ with respect to the
obvious inner product over $S^{n-1} \times \RR$.  For $\psi: S^{n-1} \times \RR
\to \RR$ the back-projection $\check{\psi}$ is
\begin{equation}	
    \check{\psi}(x) = \cR^\# \psi(x)
    =  \int_{S^{n-1}} \psi(\omega,
    \omega \cdot x) \, \textrm{d}S(\omega),
    \label{eq:backprojection}
\end{equation}
where $\textrm{d}S$ is the measure on $S^{n-1}$. 

The Radon transform $\mathcal{R}$ is a linear one-to-one map between 
$\cS(\RR^n)$ and $\cS_H(\PP^n)$ \cite{Helgason2011} in which $\cS(\RR^n)$
denotes the Schwartz class and 
\begin{equation}
\cS_H(\PP^n) = \left. 
    \begin{cases}
       &\text{for each } k \in \ZZ^+, \int_\RR F(\omega,p)p^k\, \mathrm{d}p
\text { is} \\
F \in \cS(\PP^n):   &\text{a homogeneous polynomial in } \omega_1, ... , \omega_n\\
    &\text{of degree } k 
    \end{cases}
\right\}.
\end{equation}
The correspondence can be naturally extended to distributions, and we refer the 
reader to standard references for further details.

The Radon transform has a remarkable property, that it intertwines a partial
derivative with a univariate derivative. The $i$-th partial derivative $\p / \p
x_i$ of $\varphi$ is now transformed to the derivative of $\hat{\varphi}$ with
respect to $s$ multiplied by $\omega_i$,
\begin{equation}
	\prn{\frac{\p}{\p x_i} \varphi(x)}^{\wedge} =  
	\omega_i \frac{\p}{\p s}\hat{\varphi} (\omega, s).
    \label{eq:rt}
\end{equation}
This is the key property that allows us to transform a multi-dimensional
hyperbolic problem into a collection of one-dimensional problems.  For example,
let us apply the Radon transform to the transport equation in $\RR^2$, in which
the scalar state variable $q: \RR^+ \times \RR^2 \to \RR$ satisfies,
\begin{equation}
    q_t + \theta \cdot \grad q = 0
    \quad \text{ where } \quad
    \theta \in S^1.
    \label{eq:transport2d}
\end{equation}
The transformation produces a family of 1D advection equations
\begin{equation}
    \hat{q}_t + (\theta \cdot \omega) \hat{q}_s = 0,
    \label{eq:transport1d}
\end{equation}
whose coefficient varies for each $\omega$.  Similarly, consider the
acoustic equations for $p,u,v: \RR^+ \times \RR^2 \to \RR$, where the state
variable $p$ denotes pressure, $u$  the velocity in $x_1$-direction, $v$
the velocity in $x_2$-direction,
\begin{equation}
\begin{bmatrix} p \\ u \\ v \end{bmatrix}_t 
+ \begin{bmatrix} 0 & K_0 & 0 \\ 1/\rho_0 & 0 & 0 \\ 0 & 0 & 0 \end{bmatrix} 
\begin{bmatrix} p \\ u \\ v \end{bmatrix}_{x_1}
+ \begin{bmatrix} 0 & 0 & K_0 \\ 0 & 0 & 0 \\ 1/\rho_0 & 0 & 0 \end{bmatrix} 
\begin{bmatrix} p \\ u \\ v \end{bmatrix}_{x_2}
    = 0.
    \label{eq:acoustic2d}
\end{equation}
After the transform, we obtain 
\begin{equation}
\begin{bmatrix} \hat{p} \\ \hat{u} \\ \hat{v} \end{bmatrix}_t
+ \begin{bmatrix} 0 & \omega_1 K_0  & \omega_2 K_0  \\ \omega_1 / \rho_0 & 0 & 0
\\ \omega_2/ \rho_0 & 0 & 0 \end{bmatrix} \begin{bmatrix} \hat{p} \\ \hat{u} \\
\hat{v} \end{bmatrix}_s = 0.
\end{equation}
This PDE has one spatial dimension in variable $s$. Letting $\mu = \omega_1 u +
\omega_2 v$  and  $\nu = -\omega_2 u + \omega_1 v$, \eqref{eq:acoustic2d} can
be rewritten as three equations for new states $\hat{p},\hat{\mu}$ and
$\hat{\nu}$. If one omits the trivial equation $\nu_t = 0$, the equation
\eqref{eq:acoustic2d} is reduced to the 1D acoustic equations,
\begin{equation}
    \begin{bmatrix} 
        \hat{p} \\ \hat{\mu}
    \end{bmatrix}_t
    + 
    \begin{bmatrix}
    0 & K_0 \\ 1/\rho_0 & 0
    \end{bmatrix}
    \begin{bmatrix}
        \hat{p} \\ \hat{\mu}
    \end{bmatrix}_s
    = 0. 
    \label{eq:acoustic1d}
\end{equation}
In this case, the equation depends on $\omega$ through the variable $\mu$.
However, the equation itself is invariant over all $\omega$, owing to the fact
that  the problem \eqref{eq:acoustic2d} is isotropic.  Moreover, note that this
is exactly the same equation obtained in the physical space if you consider the
case of a plane wave where the data varies only in the direction $\omega$ so
that derivatives in the orthogonal direction vanish.

The Radon transform therefore transforms $n$-dimensional hyperbolic problems
such as \eqref{eq:transport2d} and \eqref{eq:acoustic2d} into their
$1$-dimensional counterparts \eqref{eq:transport1d} and \eqref{eq:acoustic1d},
respectively.

\subsection{Multi-dimensional extension of large time-step (LTS) methods}
\label{sec:linearlts}

Previous dimensional splitting methods
\cite{acta_splitting,strangsplit,Crandall1980,fvmbook} such as Strang splitting
do not allow a natural extension of large time-step (LTS) methods
\cite{largetimestep1, largetimestep2, largetimestep3} to multiple spatial
dimensions. In order to take large time-steps for constant coefficient
multi-dimensional hyperbolic problems, one can use the Fourier transform, for
example.  Upon taking the Fourier transform, one is left with a set of ordinary
differential equations (ODEs) different from the original problem
\cite{spectral-trefethen, spectral}. On the other hand, using the Radon
transform, one obtains a dimensional splitting that reduces the
multi-dimensional problem into a family of one-dimensional counterparts of
similar (if not identical) form.  This allows 1D LTS methods to be applied for
each of these problems, and the multi-dimensional solution is obtained by
computing the inverse of the Radon transform.  Moreover, the Radon transform
provides an intuitive geometrical interpretation as a decomposition into planar
waves and thus yields other useful applications.  These applications will be
illustrated in Section \ref{sec:radonsplituse}.

This multi-dimensional extension of the LTS method for the constant coefficient
case is very straightfoward. Taking the Radon tranform of the problem as above,
one obtains a set of $1$D problems such as \eqref{eq:transport1d} or
\eqref{eq:acoustic1d}. Then one applies the $1$D LTS solution operator $\cK$ to
evolve the initial data $\hat{u}_0(\omega,s)$ for each $\omega$ up to
desired final time $T$. The operator $\cK$ may depend on the direction
$\omega$, so we denote the dependence as a parameter by writing $\cK =
\cK(T;\omega)$. This yields the Radon transform of the solution at time
$T$,
\begin{equation}
    \hat{q}(T,\omega,s) = \cK(T;\omega) \hat{q}_0(\omega,s).
    \label{eq:solveRu}
\end{equation}
Then, to compute the solution $q$ we can apply the inversion formula
\begin{equation}
    c_n q(T,x) = 
    \left\{
    \begin{aligned}
    &\cR^\# \frac{d^{n-1}}{ds^{n-1}} \, \hat{q}(T,\omega,s)  
                    & \text{ if } n \text{ is odd,}\\
    & \cR^\# H_s \frac{d^{n-1}}{ds^{n-1}} \, \hat{q}(T,\omega,s) 
                    & \text{ if } n \text{ even,}\\
    \end{aligned}
    \right.
    \label{eq:inverse}
\end{equation}
where the constant $c_n = (4 \pi)^{(n-1)/2} \Gamma(n/2)/\Gamma(1/2)$ and $H_s$
denotes the Hilbert transform.  Much is known about the inversion; see standard
texts such as \cite{Helgason2011,ctbook} for more details.

This splitting can also be related to the Strang splitting, if one views it as
a decomposition of the multi-dimensional problem into planar wave propagation.
In Strang splitting one constructs the planar waves emanating in varying
directions by dividing a single time-step into multiple successive planar wave
propagations.  The Radon transform decomposes the multi-dimensional directions
by explicitly discretizing the sphere $S^{n-1}$.

\begin{figure}
    \centering
    \begin{tabular}{m{0.25cm} m{0.35\textwidth} m{0.45\textwidth}}
        \rotatebox{90}{$t = 0$} &
        \includegraphics[width=0.35\textwidth]{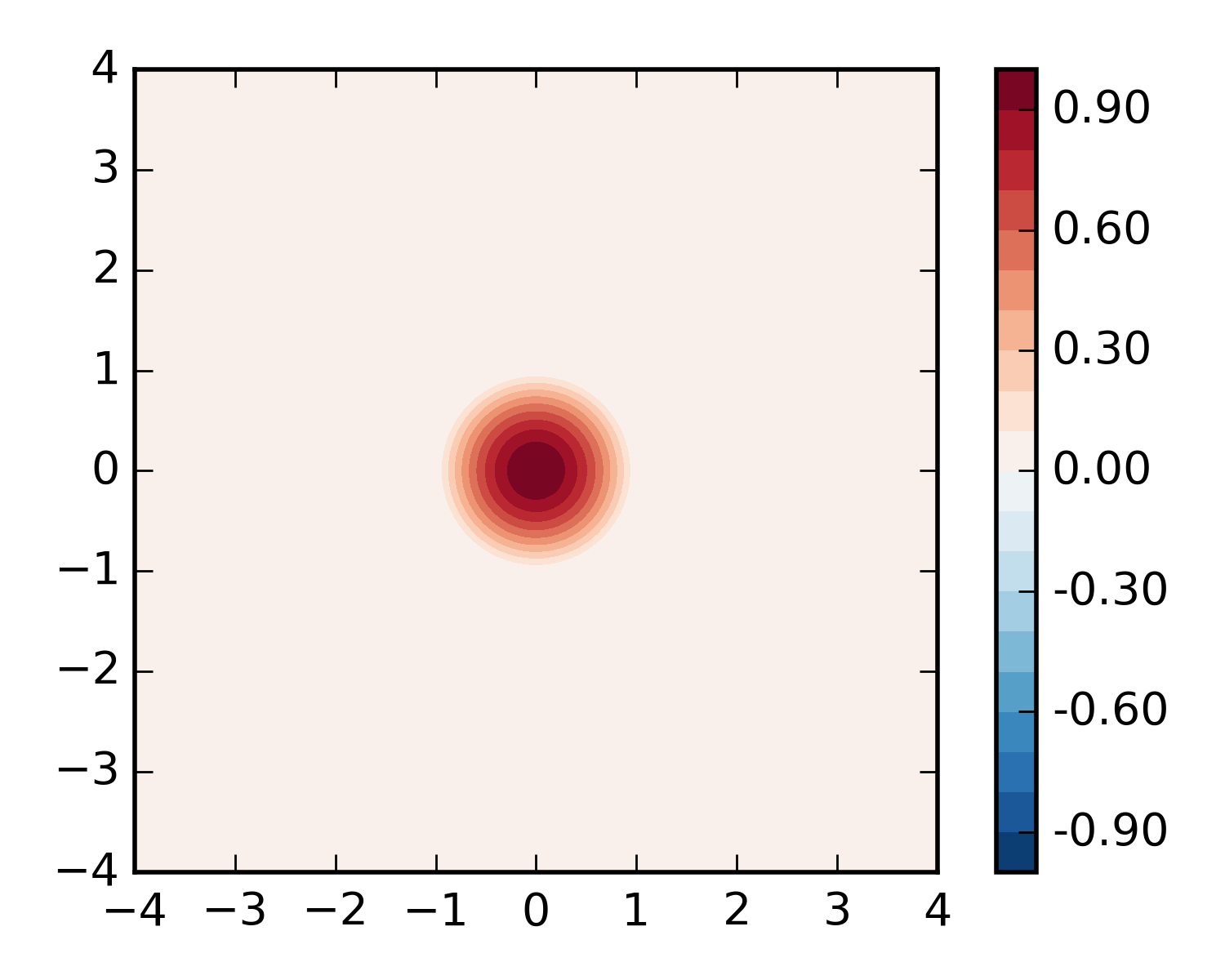}&
        \vspace{1.5em}
        \includegraphics[width=0.45\textwidth]{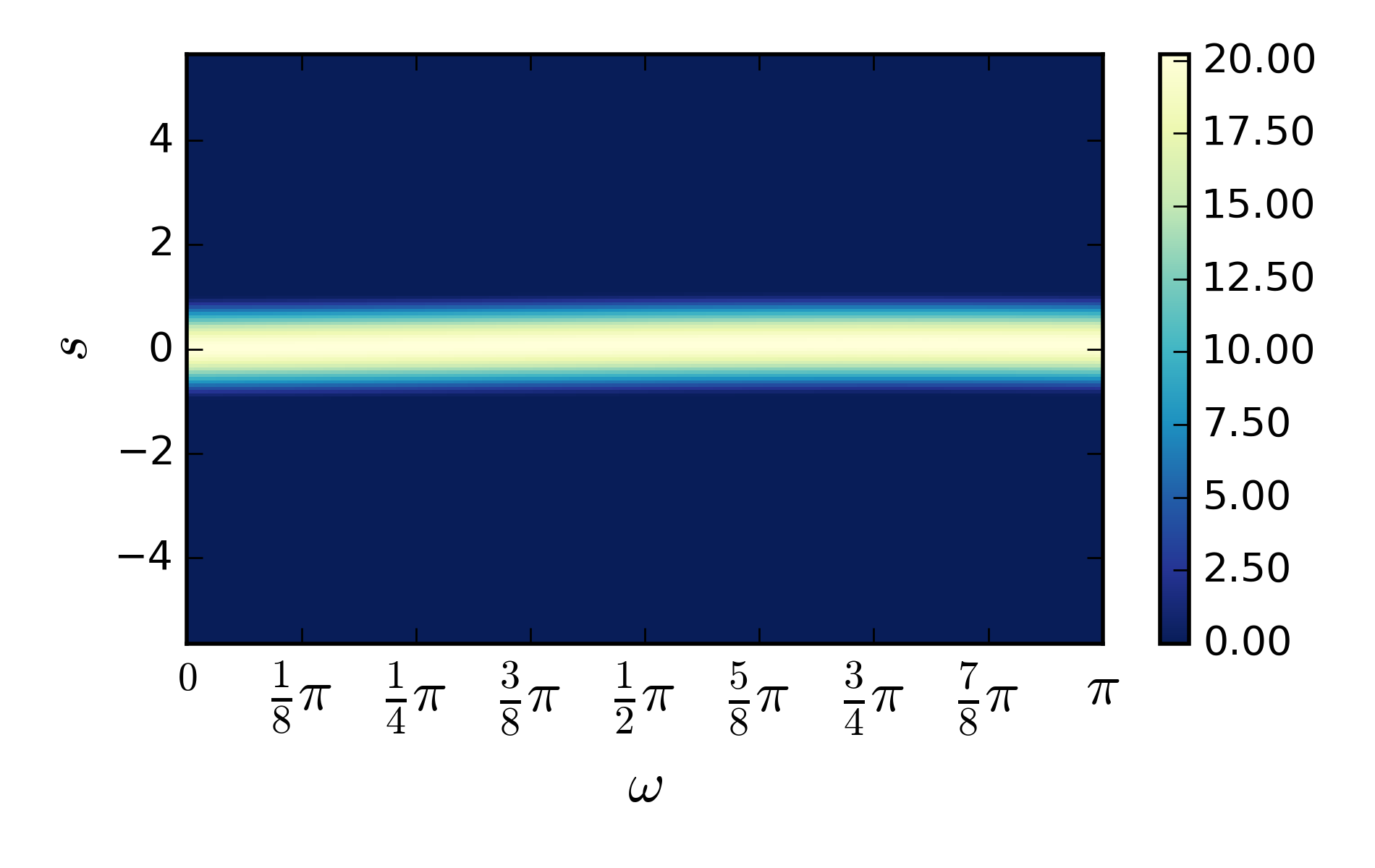}\\
        \rotatebox{90}{$t = 1$} &
        \includegraphics[width=0.35\textwidth]{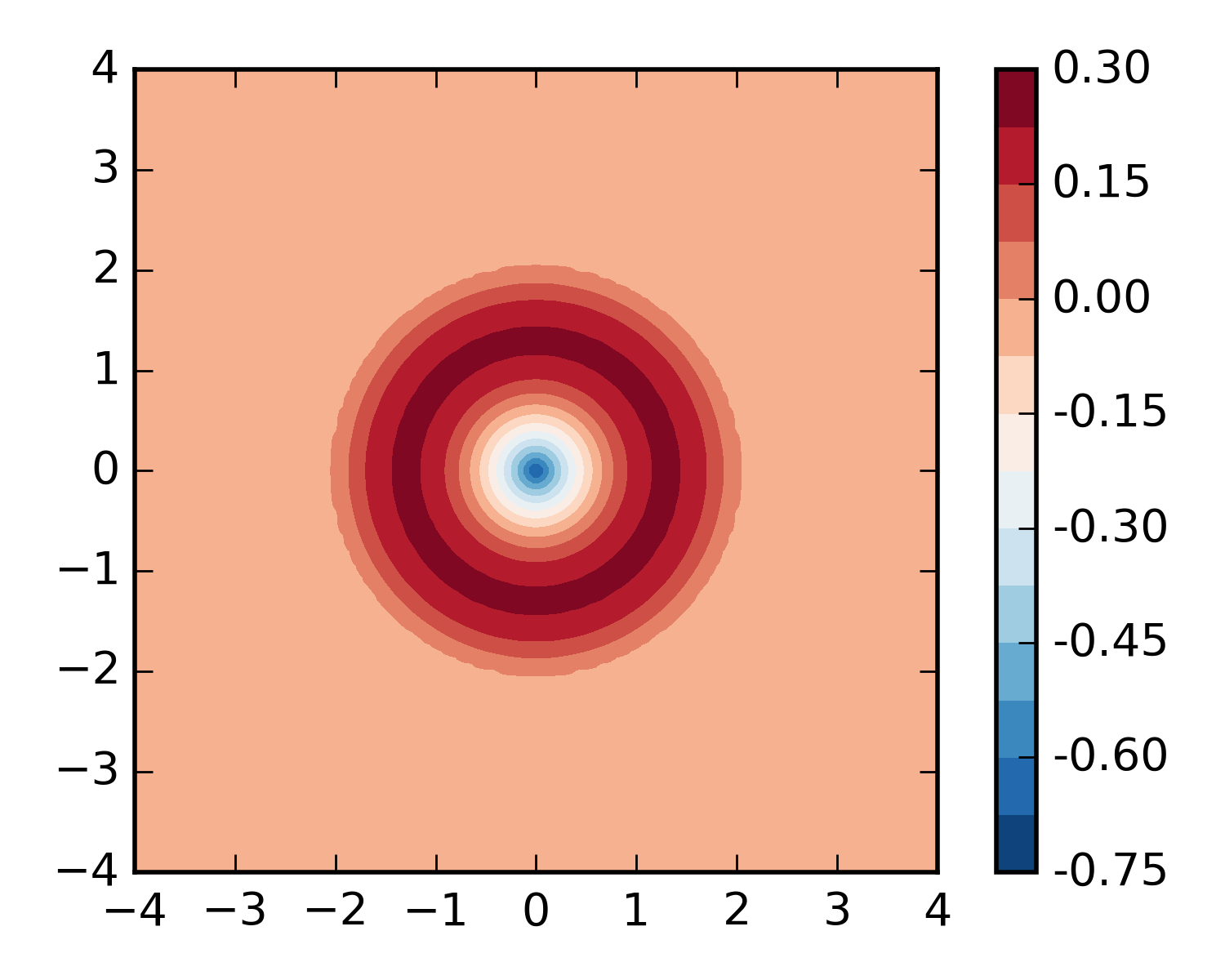}&
        \vspace{1.5em}
        \includegraphics[width=0.45\textwidth]{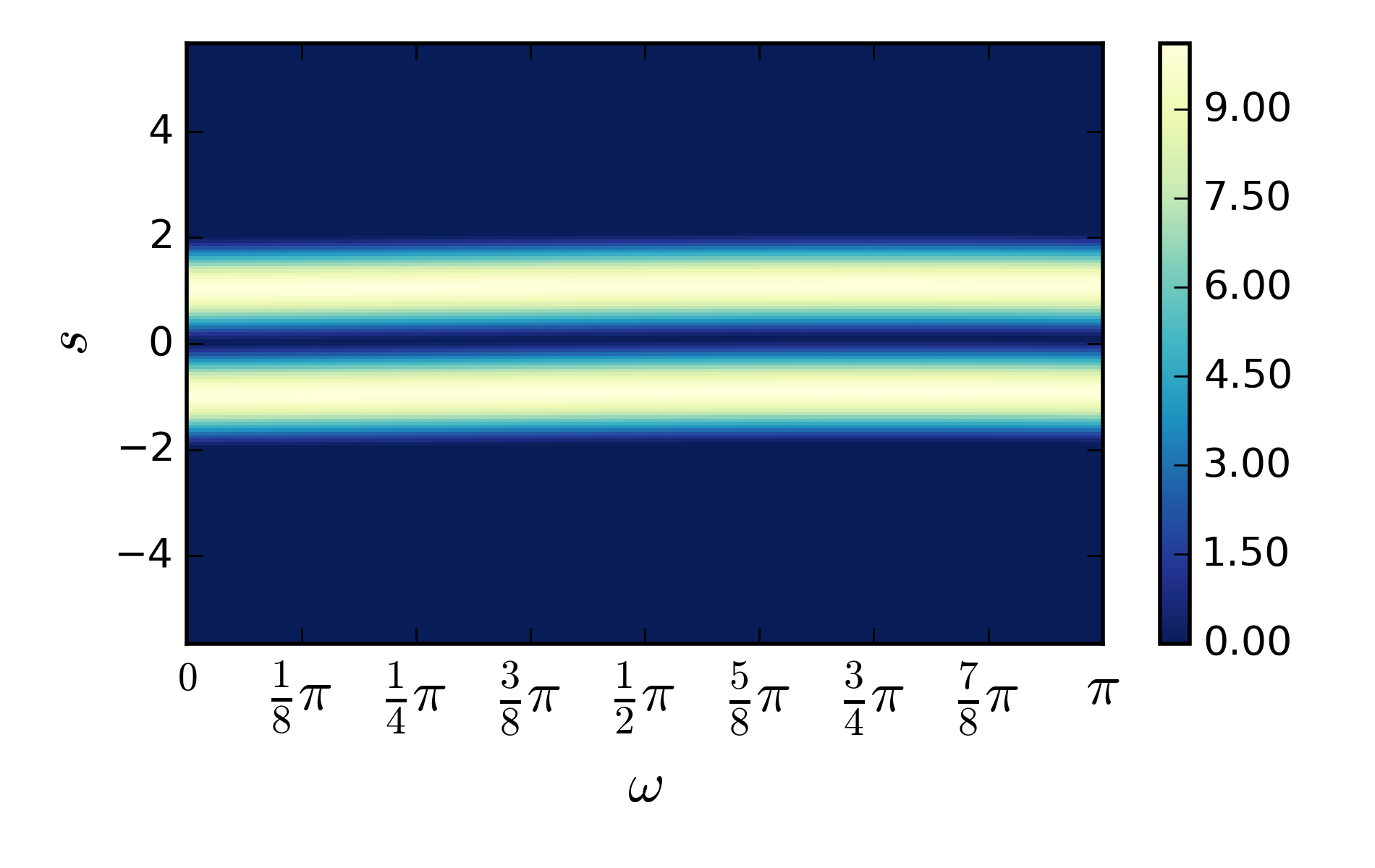}\\
        \rotatebox{90}{$t = 3$} &
        \includegraphics[width=0.35\textwidth]{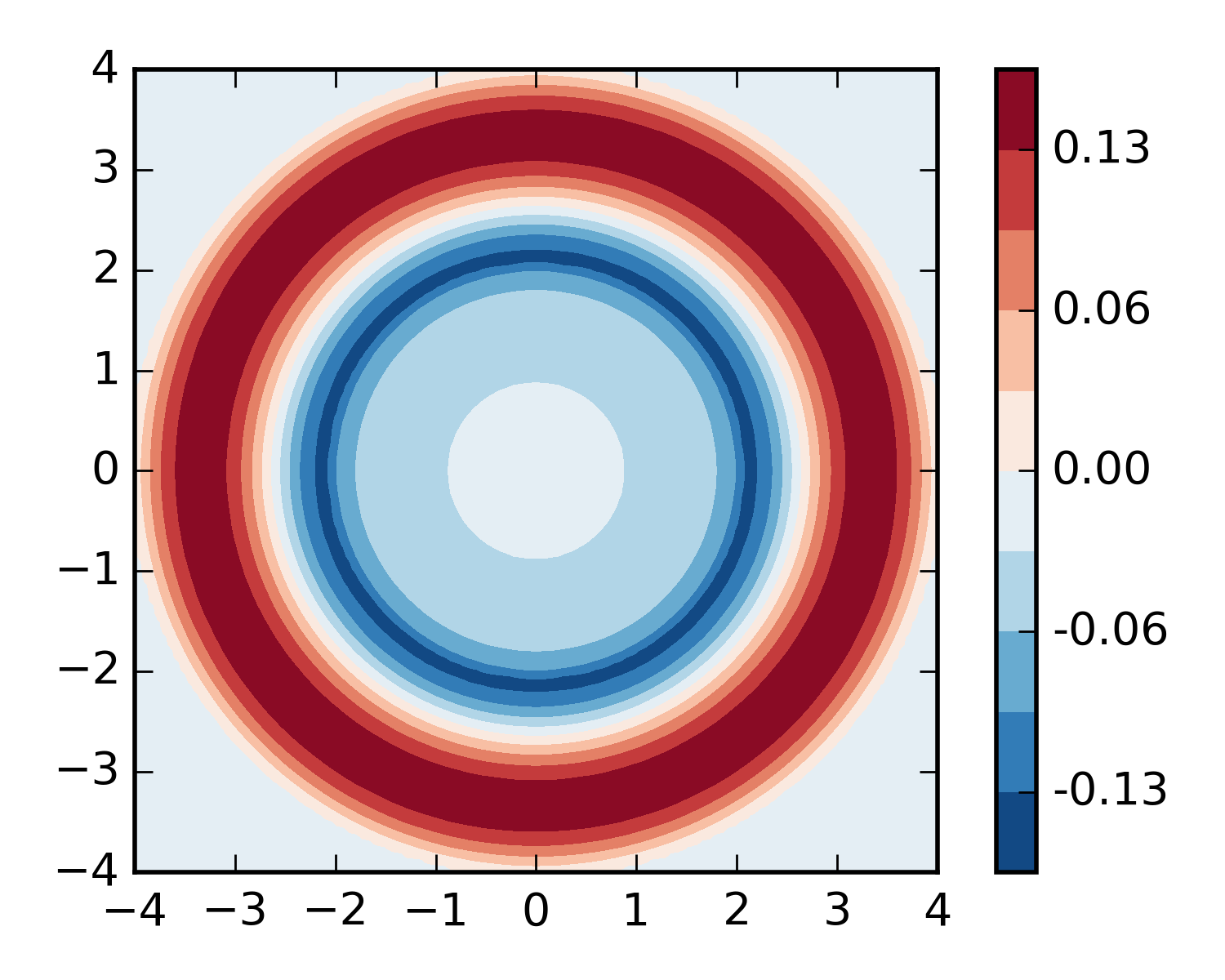}&
        \vspace{1.5em}
        \includegraphics[width=0.45\textwidth]{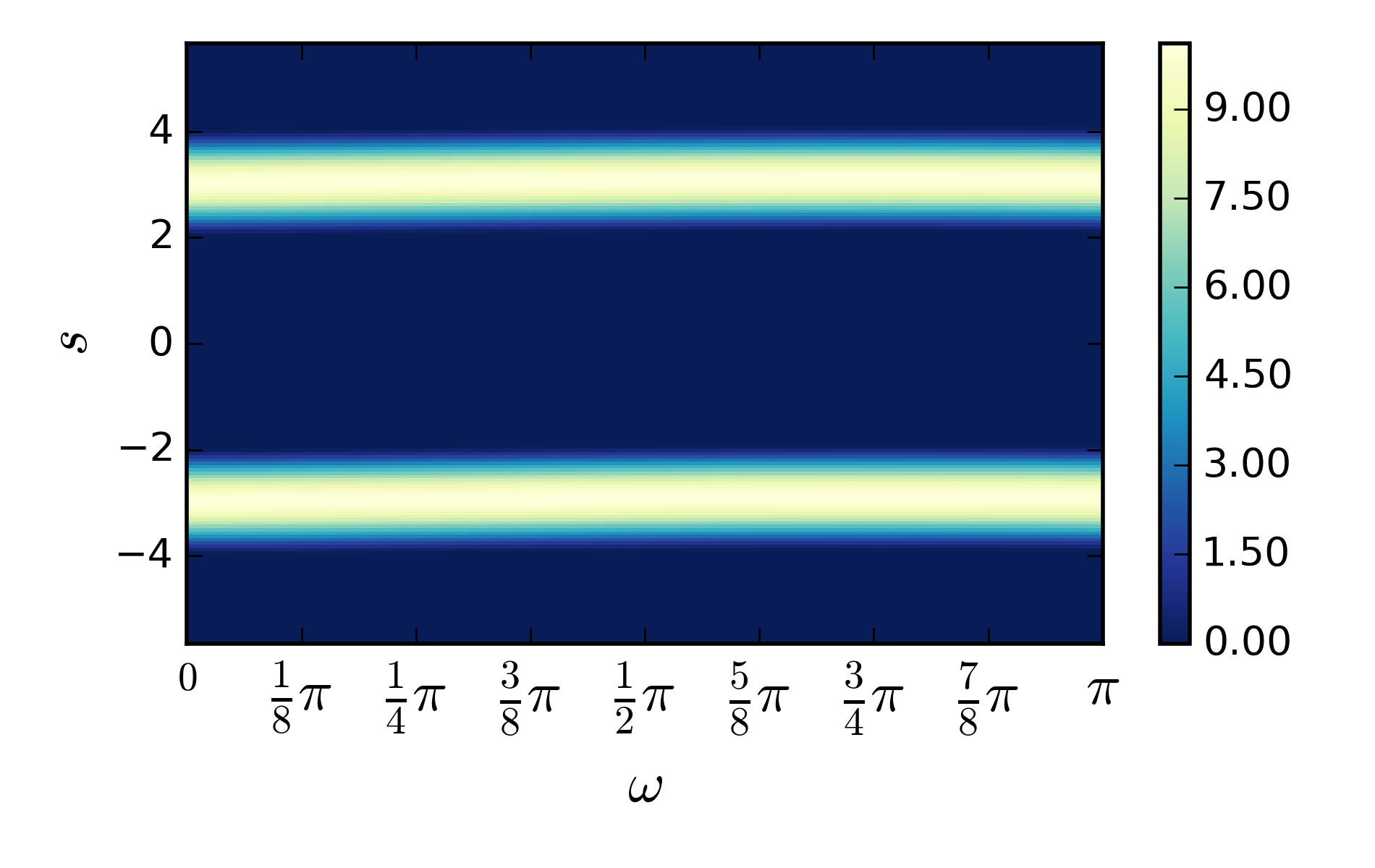}\\
    \end{tabular}
    \caption{The solution to the acoustic equation using the Radon transform
        in the square domain $[-4,4] \times [-4,4]$.
    The pressure $p$ is shown on the left column and its continuous
    Radon transform $\hat{p}$ is shown on the right column, 
    at times $t=0$ (first row), $t=1$ (second row), 
    and $t=3$ (third row).
    }
    \label{fig:acoustic}
\end{figure}

Let us consider a concrete example, the 2D acoustic equation
\eqref{eq:acoustic2d}.  Let us set $K_0 = \rho_0 = 1$, so that we have the
sound speed $c = 1$, and impose the initial conditions
\begin{equation}
    \begin{aligned}
    q_0(x) &= \begin{bmatrix} 
                   p_0(x_1,x_2) \\
                   0  \\
                   0  \\
               \end{bmatrix}, \\
    p_0(x) &=
    \left\{
    \begin{aligned}
      \cos(\pi(x_1^2 + x_2^2)/2) 
                      &\quad \text{ if } \, x_1^2 + x_2^2 < 1, \\
      0               &\quad \text{ otherwise. }  \\
    \end{aligned}
    \right.
    \end{aligned}
    \label{eq:coshump}
\end{equation}
The initial pressure profile is a cosine hump  supported in a disk of radius
$1$ centered at the origin, and the initial velocity profile is identically
zero. We will also set absorbing boundary conditions in the manner to be
described in Section \ref{sec:absorb}.

On the transformed side \eqref{eq:acoustic1d}, the evolution for any fixed
direction  $\omega \in S^1$ is given by the d'Alembert solution
\eqref{eq:dalembert},
\begin{equation}
    \hat{q}(t,\omega,s) = 
\frac{1}{2} \left( r_1 \hat{p}_0(\omega,s - t)  + r_2 \hat{p}_0(\omega,s + t)  \right)
\quad \text{ where } 
r_1 = \begin{bmatrix} 1 \\ 
               \omega_1 \\ 
               \omega_2 \end{bmatrix} \text{ and }
r_2 = \begin{bmatrix} 1 \\
              -\omega_1 \\ 
              -\omega_2 \end{bmatrix}.
    \label{eq:dalembert}
\end{equation}
This reduces to simple shifts at corresponding speeds, which can be computed
easily up to any time $t$. This is precisely the 1D LTS solution for the
constant coefficient case.

The solution to the acoustic equation computed on the domain $[-4,4] \times
[-4,4]$ is shown in the left column of Figure \ref{fig:acoustic}.  The Radon
transform of the pressure term $\hat{p}$ is plotted in the right column of the
same figure. Note that this problem is radially symmetric about the origin. A
consequence of this is that the Radon transform is invariant with respect to
the variable $\omega$, hence the Radon transform of the solutions at different
times all appear as horizontal stripes.  (There is a small amount of shift,
following from the fact that for an image of even size $N$, the origin is
chosen as the $(N/2,N/2)$-pixel, slightly off center.)

A key observation is that the evolution of the solution in the transformed
variables is a sum of two shifting horizontal stripes, although the wave
profile in the spatial domain propagates radially.  For each fixed angle
$\omega$, one only need solve the d'Alembert solution \eqref{eq:dalembert},
which is easy to solve to any time $t$ by shifting the initial profile twice
each according to two opposite speeds, and summing them. Intuitively, the
shifts correspond to the propagation of decomposed planar waves for any fixed
normal directions in $S^1$.

The actual computational did not make use of the continuous Radon transform
\eqref{eq:radon}, but rather a completely discrete approximation called the
DRT, which will be introduced and discussed in further detail in Section
\ref{sec:drt}.  Here it will suffice to mention that a grid of size $128 \times
128$ was used and prologation of $p=2$ was used for the DRT, and that the
continuous transform can be obtained by an easy change of variables
\eqref{eq:drt2rt} which scale the domain and amplitude of the DRT, and that the
change of variables do not affect the intertwining property. The 1D LTS method
can still be used on the DRT, just as in the case of the continuous transform.
The computational cost for this solution is conjectured to be $\cO(N^{5/2} \log
N)$: $\cO(N^2 \log N)$ for the forward DRT, $\cO(N^2)$ for the 1D LTS solution,
and $\cO(N^{5/2} \log N)$ for the inverse DRT (see Section \ref{sec:drt}). We
note that the 1D LTS method can be applied for each angle in parallel.

\begin{figure}
    \centering
    \begin{tabular}{m{0.25cm} m{0.35\textwidth} m{0.45\textwidth}}
    \rotatebox{90}{$t = 0$} &
    \includegraphics[width=0.35\textwidth]{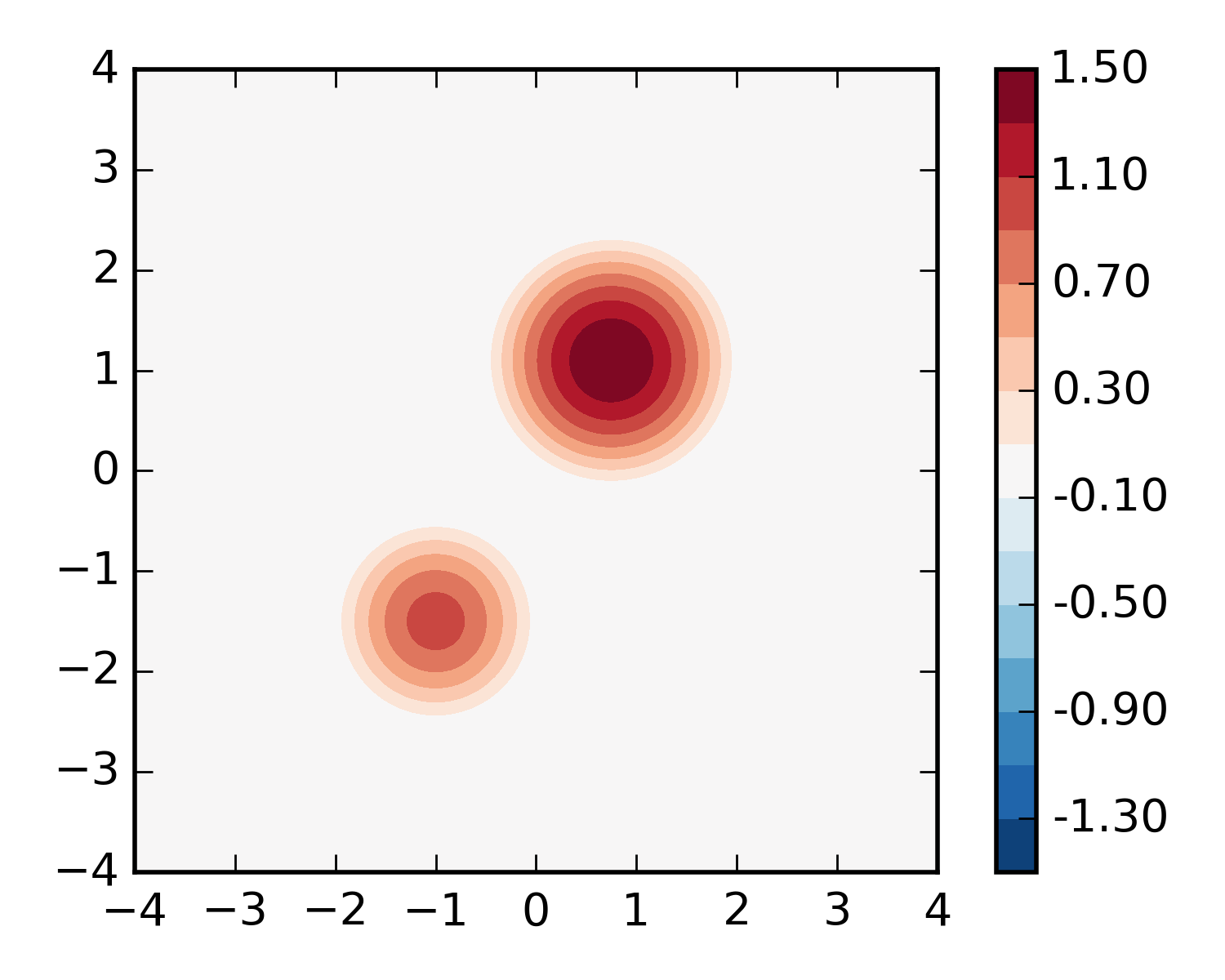}
    &
        \vspace{1.5em}
    \includegraphics[width=0.45\textwidth]{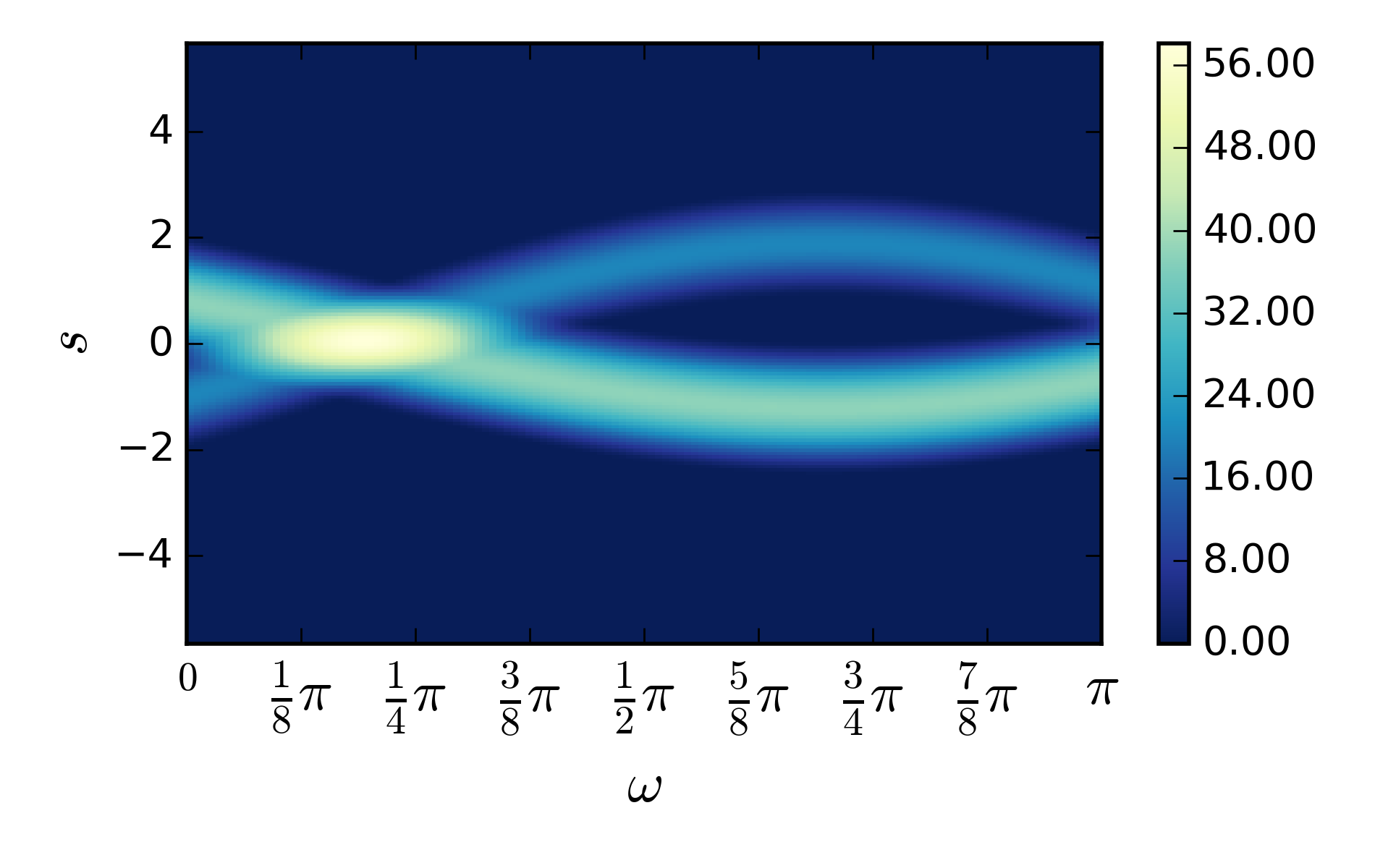}
    \\
    \rotatebox{90}{$t = 0.5$} &
    \includegraphics[width=0.35\textwidth]{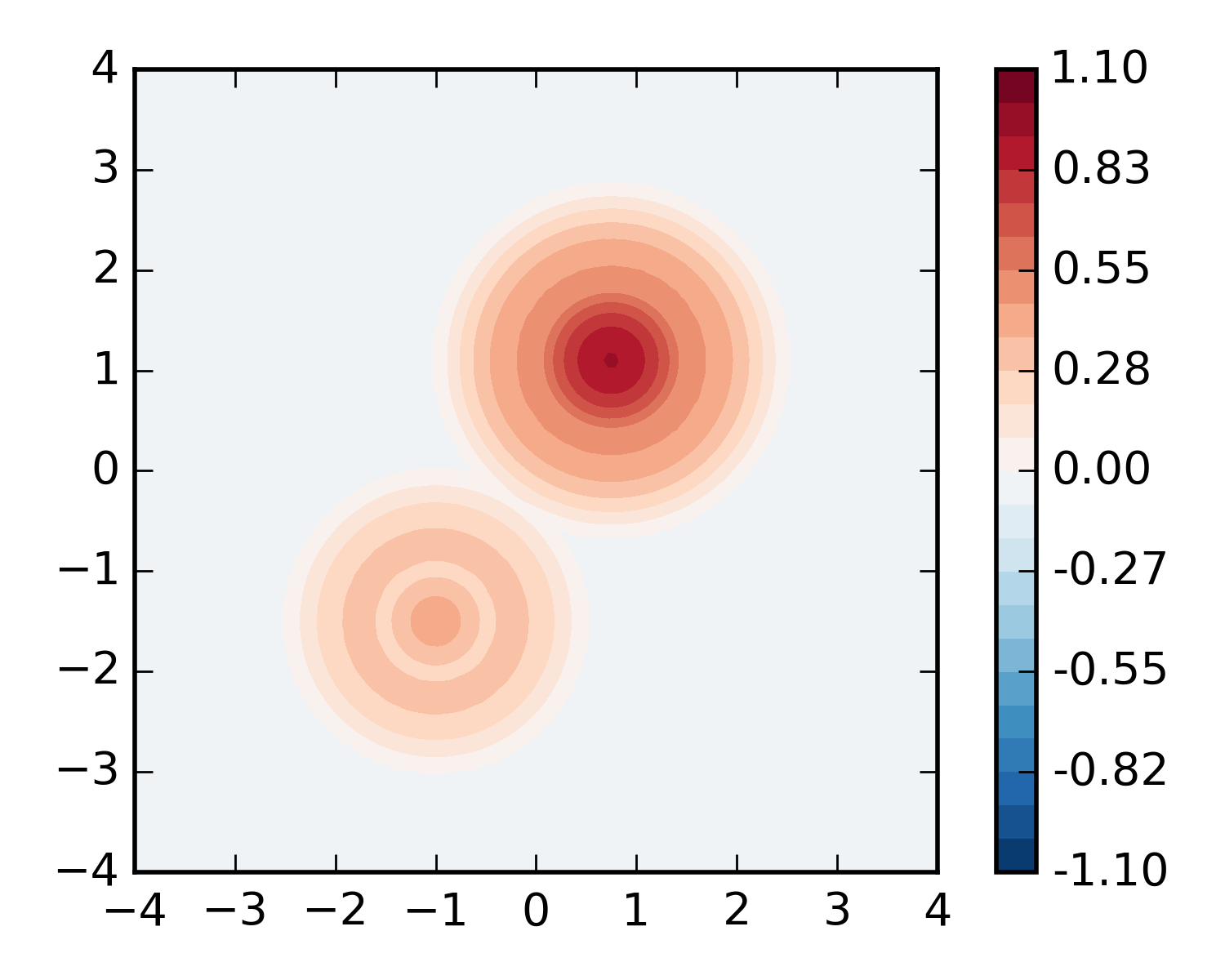}
    &
        \vspace{1.5em}
    \includegraphics[width=0.45\textwidth]{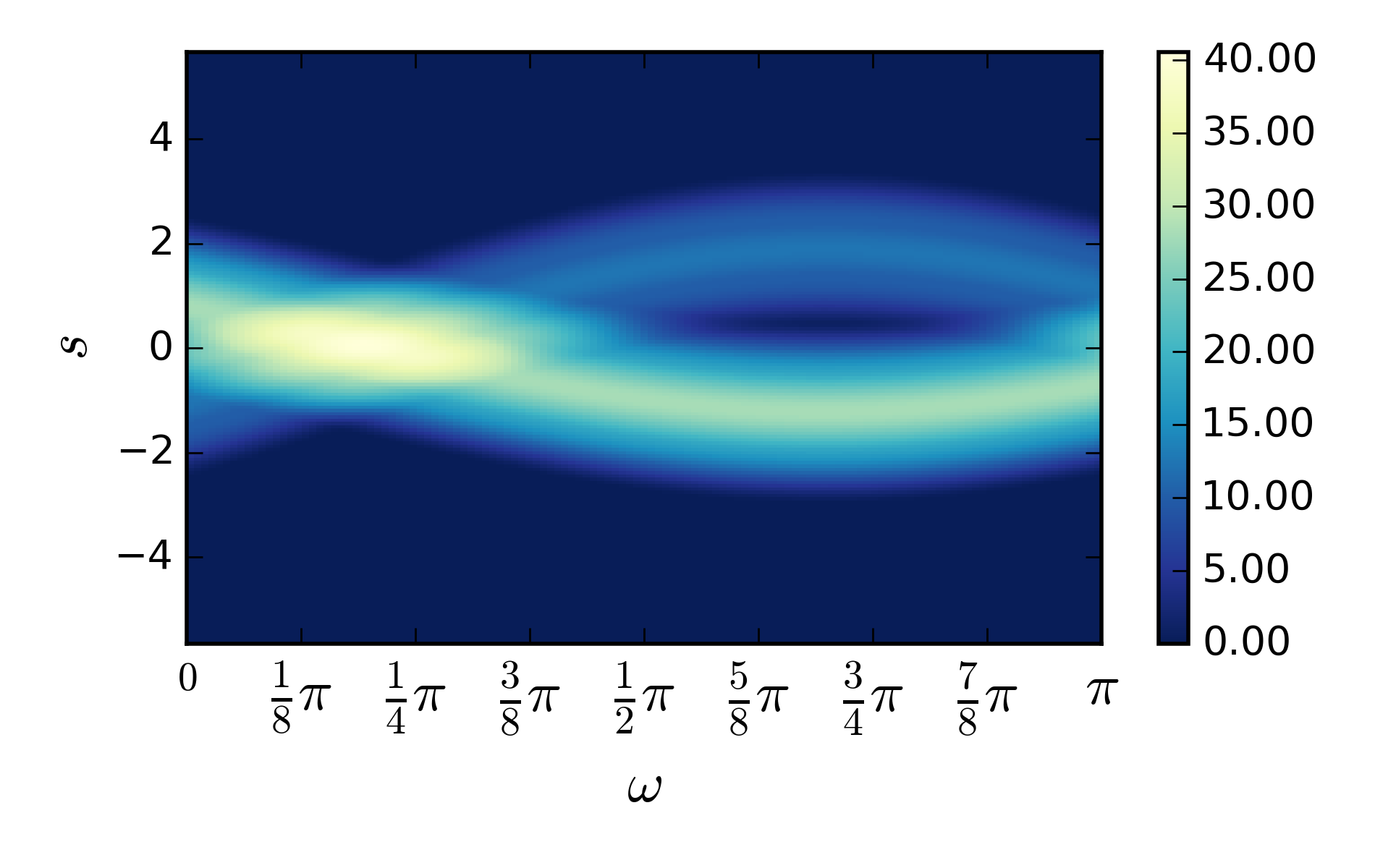}
    \\
    \rotatebox{90}{$t = 1$} &
    \includegraphics[width=0.35\textwidth]{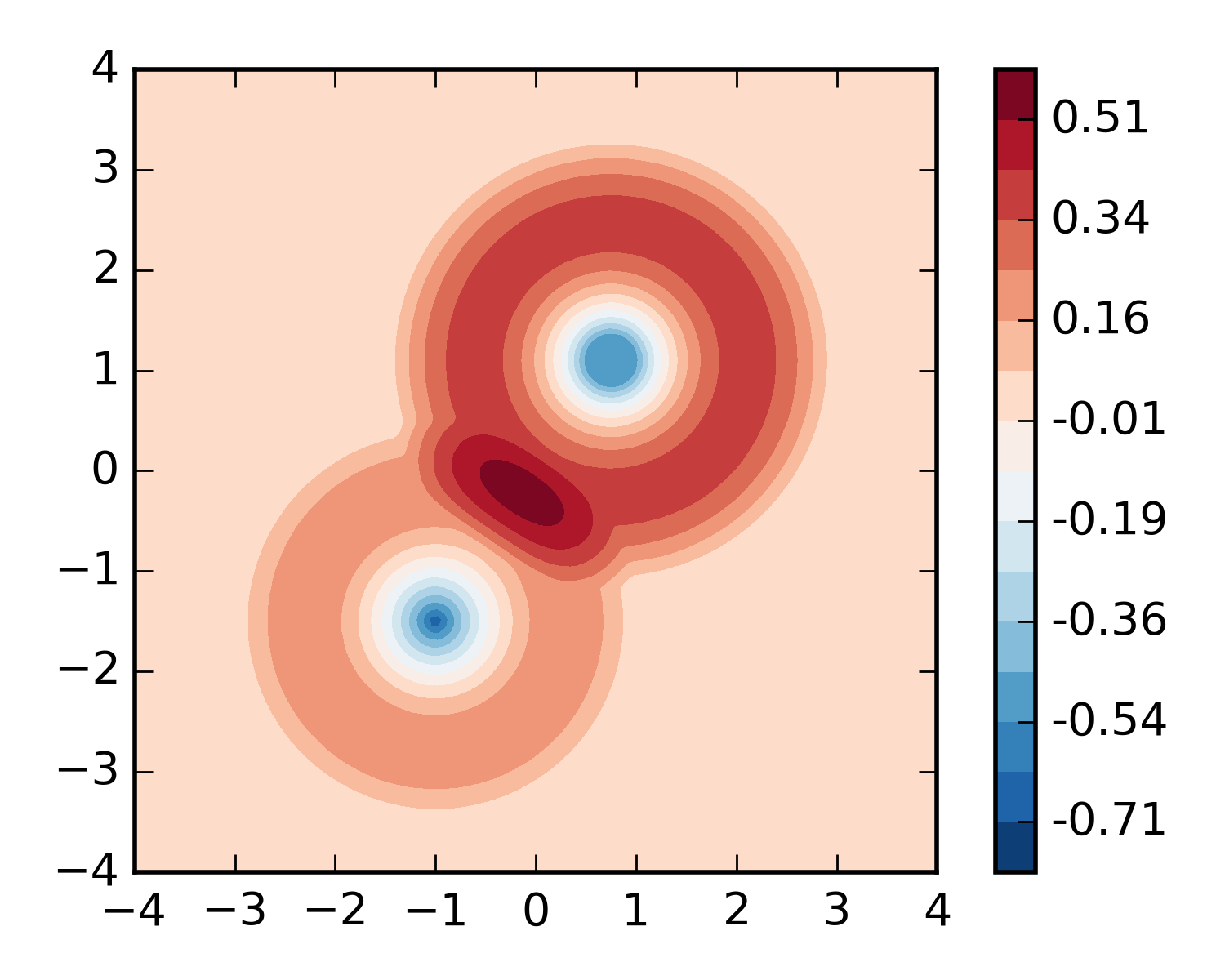}
    &
        \vspace{1.5em}
    \includegraphics[width=0.45\textwidth]{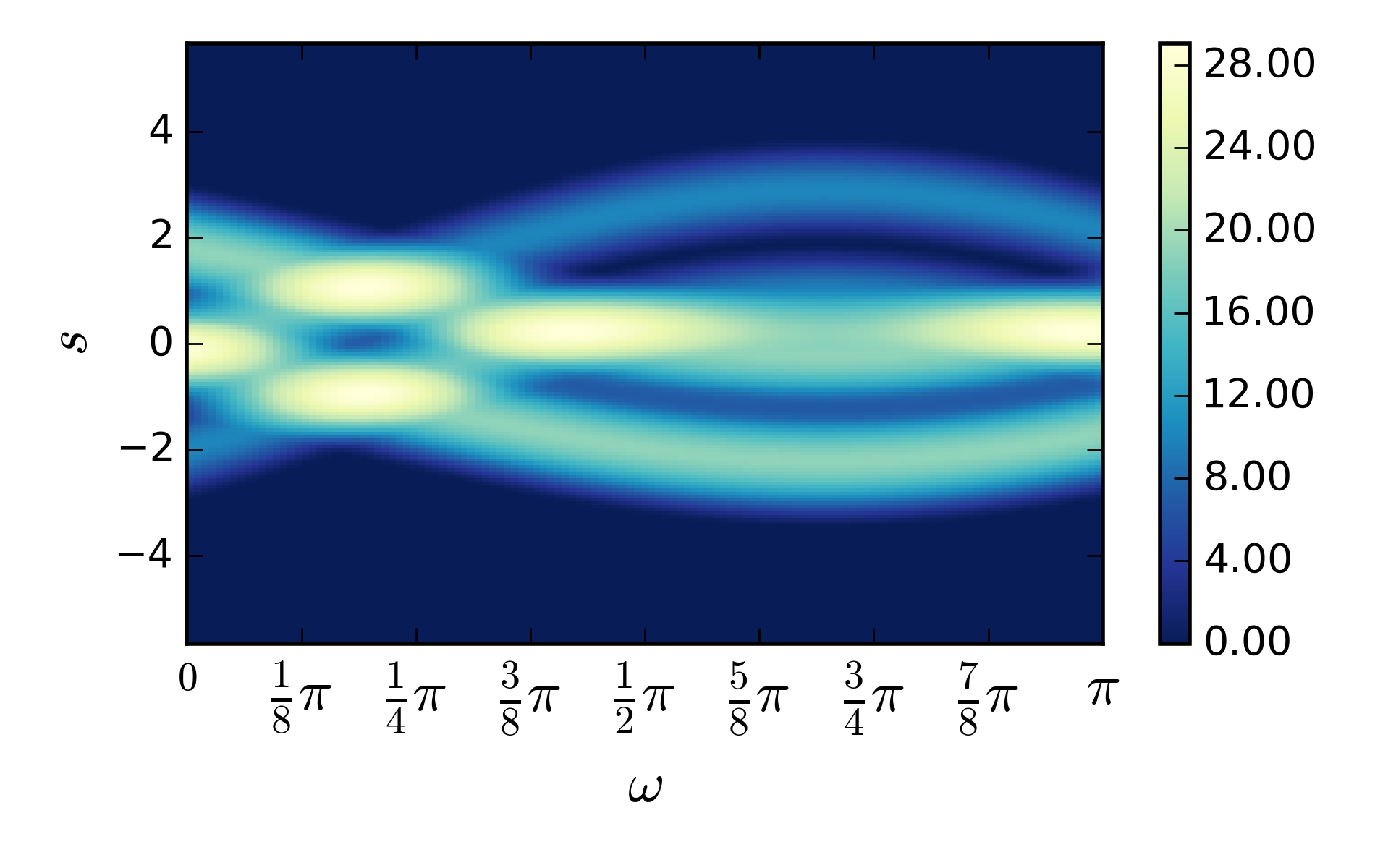}
    \\
    \rotatebox{90}{$t = 1.5$} &
    \includegraphics[width=0.35\textwidth]{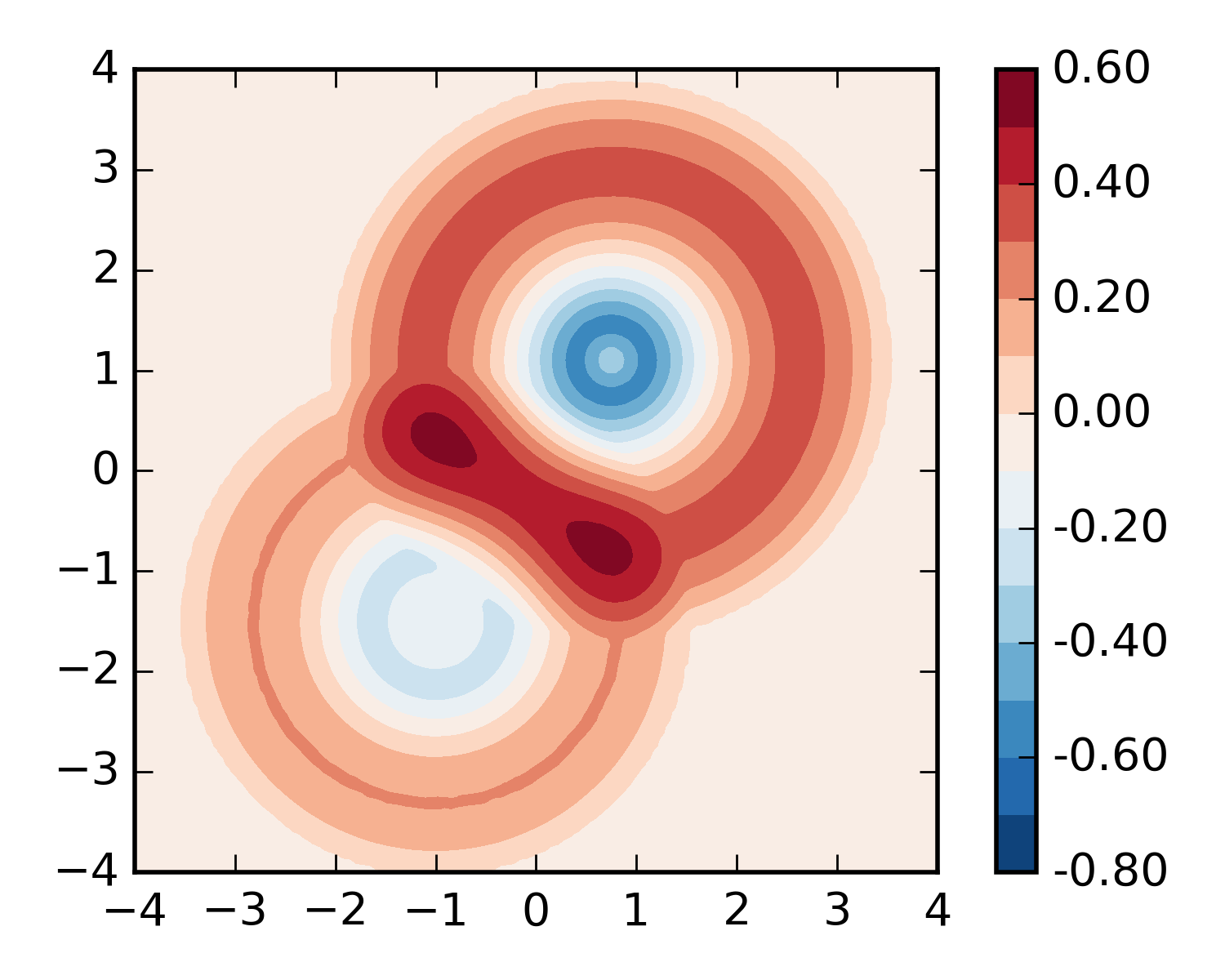}
    &
        \vspace{1.5em}
    \includegraphics[width=0.45\textwidth]{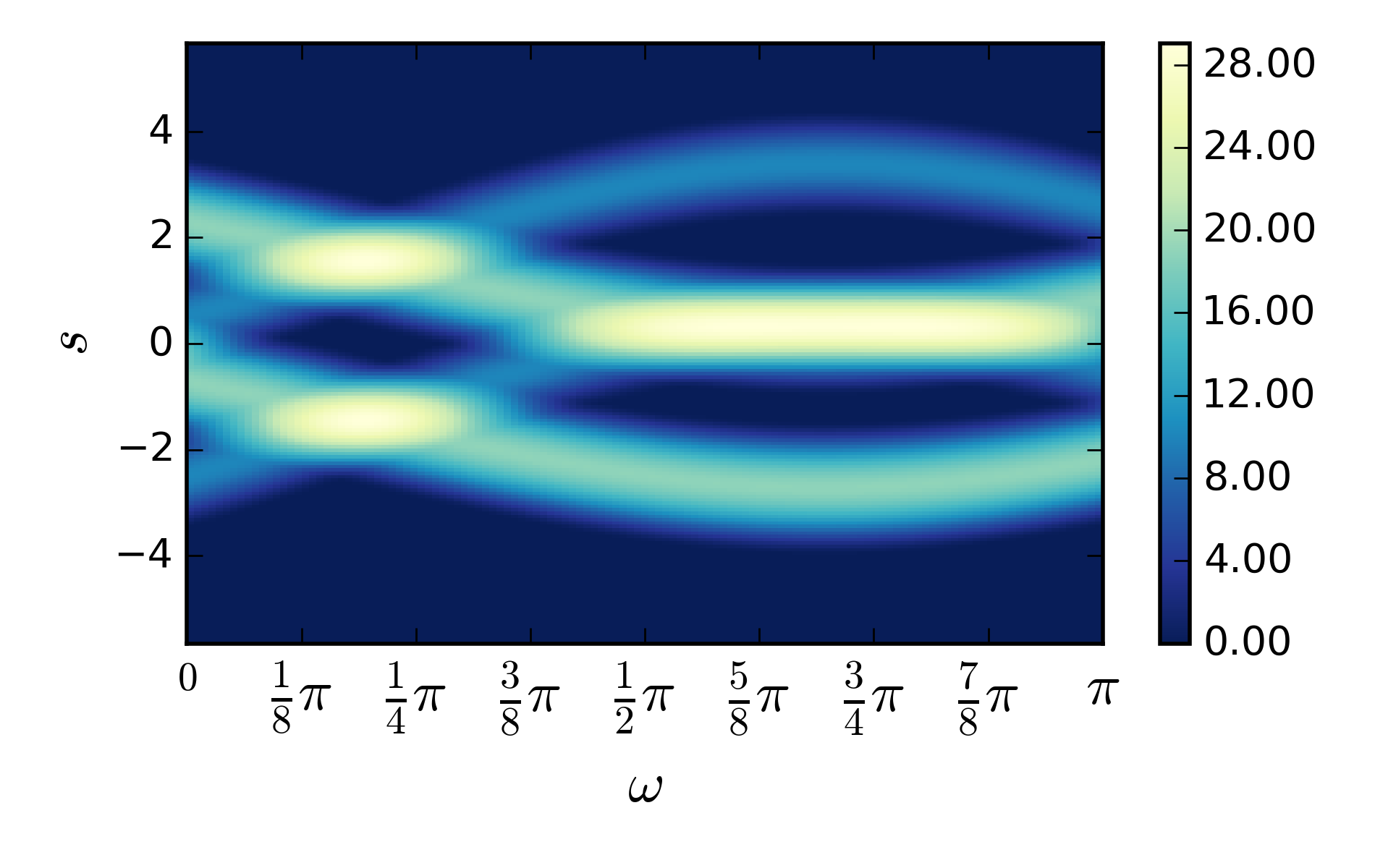}
    \\
    \end{tabular}
    \caption{The solution to the acoustic equation using the Radon transform
    in the square domain $[-4,4] \times [-4,4]$.
    The pressure $p$ is shown on the left column and its continuous
    Radon transform $\hat{p}$ is shown on the right column, 
    at times $t=0$ (first row), $t=0.5$ (second row), $t=1$ (third row),
    and $t=1.5$ (fourth row).}
    \label{fig:acoustic_2humps}
\end{figure}

\begin{figure}
    \centering
    \begin{tabular}{ c}
    \includegraphics[width=0.7\textwidth]{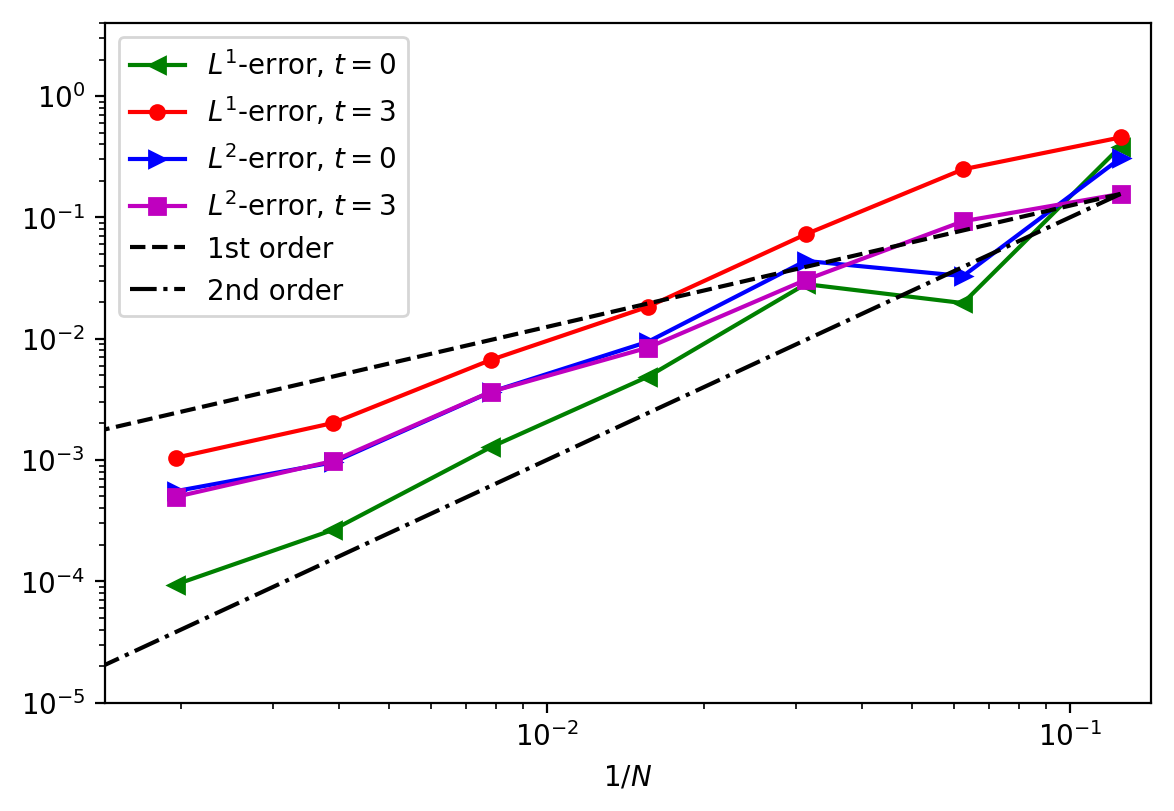}
    \\
    \begin{tabular}{rrrr}
    \hline
    \multicolumn{2}{c}{$L^1$ error }
    &
    \multicolumn{2}{c}{$L^2$ error} 
    \\
    \hline
        $t=0$ & $t=3$  &  $t=0$  &  $t=3$\\
\hline
0.38371085 &  0.45606217 & 0.30756990  & 0.15607800 \\
0.01964332 &  0.24903830 & 0.03281770  & 0.09287436 \\
0.02817384 &  0.07243586 & 0.04375375  & 0.03049178 \\
0.00486719 &  0.01834732 & 0.00944203  & 0.00843966 \\
0.00127273 &  0.00666983 & 0.00363240  & 0.00364606 \\
0.00026550 &  0.00201824 & 0.00095616  & 0.00098109 \\
0.00009389 &  0.00104277 & 0.00055271  & 0.00049592 \\
\hline
    \end{tabular}
    \end{tabular}
    \caption{Convergence plot of the pressure variable $p$ of the splitting
solution to the acoustics equation \eqref{eq:acoustic2d} with initial
conditions \eqref{eq:coshump} at times $t=0$ and $t=3$, for grid cell sizes
$N=8,16,32,64,128,256,512$ (left). The reference solution was computed by
solving an equivalent 1D problem with 4000 grid cells and the weighted
$L^1,L^2$ norms \eqref{eq:weightL2} were used to compute the error.
The numerical values are also displayed (right).}
    \label{fig:acoustic_conv}
\end{figure}

Since the problem is radially symmetric, we can compare the solution to a
1D reference solution of high accuracy. We computed the 1D problem using
Godunov flux with 4000 grid cells, as implemented in the \textsc{Clawpack}
software package \cite{clawpack}.  We compared the diagonal slice of our 2D
solution at angle $\pi/4$ with the reference solution at time $t=0$ and $t=3$.
To observe the accuracy of the solution with respect to the grid-size,
solutions of sizes $N = 8,16,32,64,128,256,512$ are also compared. The
error was computed for the pressure variable $p$ using the weighted $L^1$
and $L^2$ norms, 
\begin{equation}
\begin{aligned}
&\left(\int_0^{4\sqrt{2}} | p(\rho,t)- p_\textrm{ref}(\rho,t)|^{\mathfrak{p}}
\rho \, \textrm{d}\rho \right)^{\frac{1}{\mathfrak{p}}}\\
&\quad \quad 
\text{ where } \rho = \sqrt{x_1^2 + x_2^2}, 
\quad
p(\rho,t) = p\left( \frac{\rho}{\sqrt{2}},\frac{\rho}{\sqrt{2}},t \right),
\quad
\mathfrak{p} = 1,2.
\end{aligned}
\label{eq:weightL2}
\end{equation}
The comparison results are displayed in Figure \ref{fig:acoustic_conv}. The
error at the later time is at the level of the initial discretization error. We
also observe that the convergence rate is between first and second-order with
respect to the cell diameter $1/N$.

While this particular problem was radially symmetric, the splitting is by no
means restricted to problems with radial symmetry. Let us modify the initial
conditions above so that it is the sum of two cosine humps of different radii
and heights,
\begin{equation}
    \begin{aligned}
    q_0(x) &= \begin{bmatrix} 
                   p_1(x_1,x_2) \\
                   0  \\
                   0  \\
               \end{bmatrix},\\
    p_1(x) 
       &=            p_0(x_1 + 1 , x_2+ 1.5) 
            + 1.5 \, p_0(1.25(x_1 - 0.75),1.25(x_2 - 1.1)).
    \end{aligned}
    \label{eq:2humps}
\end{equation}
The splitting solution and its continuous Radon transform is plotted in Figure
\ref{fig:acoustic_2humps}. In the first row of the figure, the initial
condition and its Radon transform are shown. The two cosine humps in the
initial condition each correspond to a sinusoidal signal on the transformed
side.  Recall the horizontal line centered at $s=0$ from the previous example
(Figure \ref{fig:acoustic}). The sinusoidal shift away from $s=0$ is due to the
fact that translation  is an anisotropic operation. This can also be deduced
from the transformed transport equation \eqref{eq:transport1d} in which the
transport speed is $\theta \cdot \omega$, that is, $\cos \phi$ where
$\phi$ is the angle between transport direction $\theta$ and the direction
of the transform $\omega$.  For example, when the cosine hump at the origin
$p_0$ \eqref{eq:coshump} is transported away from the origin by $r(\theta
\cdot \omega)$, $\hat{p_0}$ is shifted by $\hat{p_0}(\omega, s -
r\cos(\phi))$.

In any case, the solution is still given by the d'Alembert solution
\eqref{eq:dalembert} and the acoustics equation can be solved exactly the
same way as before. The DRT used in the actual computations are plotted in
Figure \ref{fig:acoustic_2humps_drt}. Each corresponds to a continuous
transform in Figure \ref{fig:acoustic_2humps}.

\subsection{Splitting for the nonlinear case}\label{sec:nonlinearlts}
Here we will discuss how the splitting above to can be applied to a fully nonlinear 
system of hyperbolic equations.
For a state vector $q$, such a PDE is given in the form
\begin{equation}
    q_t + f(q)_{x_1} + g(q)_{x_2} = 0,
\end{equation}
where $f$ and $g$ are flux functions that can be nonlinear.
Taking the Radon transform as before, we obtain
\begin{equation}
    \hat{q}_t + \left[ \omega_1 f(q) + \omega_2 g(q) \right]^\wedge_s = 0.
\end{equation}
Let us define the \emph{directional flux function} as
\begin{equation}
h(q) = \omega_1 f(q) + \omega_2 g(q).
\label{eq:2dflux}
\end{equation}
Then one obtains the nonlinear $1$D equations,
\begin{equation}
    \hat{q}_t + h(q)^\wedge_s = 0.
    \label{eq:nonlinear1d}
\end{equation}
As in the acoustics equation \eqref{eq:acoustic1d}, the dependence on $\omega$
enters through the flux function $h(q)$, while the form of the equation is
invariant with respect to $\omega$.

As an example, consider the shallow water equations in 2D, in which $\rho,u,v :
\RR^+ \times \RR^2 \to \RR$ denote water height, velocity in the
$x_1$-direction and velocity in the $x_2$-direction, respectively,
\begin{equation}
    \begin{bmatrix}
        \rho \\
        \rho u \\
        \rho v \\
    \end{bmatrix}_t
    +
    \begin{bmatrix}
        \rho u \\
        \rho u^2 + \frac{1}{2} \bar{g}\rho^2 \\
        \rho uv
    \end{bmatrix}_{x_1}
    +
    \begin{bmatrix}
        \rho v \\
        \rho uv\\
        \rho v^2 + \frac{1}{2} \bar{g}\rho^2 
    \end{bmatrix}_{x_2}
    =
    0.
    \label{eq:shallow2d}
\end{equation}
Here $\bar{g}$ denotes the gravitational constant.
The Radon transform as above yields $1$D equation in the form
\eqref{eq:nonlinear1d},
in the transformed velocity variables $\mu = \omega_1 u + \omega_2 v$ and 
$\nu = -\omega_2 u + \omega_1 v$,
\begin{equation}
    \begin{bmatrix}
        \rho \\
        \rho\mu \\
        \rho\nu \\
    \end{bmatrix}^\wedge_t
    +
    \begin{bmatrix}
        \rho \mu \\
        \rho\mu^2 + \frac{1}{2} \bar{g}\rho^2 \\
        \rho\mu\nu
    \end{bmatrix}^\wedge_s = 0.
    \label{eq:shallow1d}
\end{equation}
Note that the first two equations of \eqref{eq:shallow1d} are just the 
shallow water equation in a single dimension in the normal direction of the
hyperplane, whereas the third equation is the conservation of momentum in 
the transverse direction.

We observe that the transformed equations resemble a finite volume
discretization.  Let us say that $\xi_{i,j}$ is a discretization of the
hyperplane $\{x \in \RR^n : x \cdot \omega_i = s_j \}$.  The
specific discretization for the hyperplanes can take on many different forms,
but here we will leave it in a general form.  We denote the approximation to
$\hat{q}(t_n,\omega_i,s_j)$ at time-step $t_n$ by $\hat{Q}^n_{i,j}$, 
\begin{equation}
    \hat{Q}_{i,j}^{n} \approx 
    \int_{\xi_{i,j}} q(t_n,x) 
    \, \textrm{d}m(x).
\end{equation}
For each fixed direction $i=i_0$, the collection of hyperplanes $\{\xi_{i_0,j}\}$ 
form a partition of the domain. We can consider these hyperplanes to be finite
volume cells. In the equation \eqref{eq:nonlinear1d} the flux function $h(q)$
assigns the flux between $\xi_{i_0,j}$ and $\xi_{i_0,j+1}$.  If the cell
boundary between $\xi_{i_0,j}$ and $\xi_{i_0,j+1}$ is denoted by
$\xi_{i_0,j+\half}$, we define the numerical flux $F_{i_0,j+\half}^n$ to be an
approximation to the flux at $\xi_{i_0,j+\half}$, valid from time-step $t_n$ to
$t_{n+1}$.  Then we have the finite volume update
\begin{equation}
    \hat{Q}^{n+1}_{i_0,j} = \hat{Q}^{n}_{i_0,j} 
    - \Delta t (F_{i_0,j+\half}^n - F_{i_0,j-\half}^n).
\end{equation}
Once these updates are made for all $i$, the updated $\hat{Q}^{n+1}_{i,j}$ are
combined through the inversion formula \eqref{eq:inverse} to yield the
numerical solution at time $t_{n+1}$. 

The dimensional splitting strategy would be to compute the numerical flux
$F^n_{i,j+\half}$ by solving only the 1D Riemann problems in the $x$ and $y$
directions.  Since the flux function $h(q)$ is a linear combination of normal
fluxes $f(q)$ and $g(q)$ \eqref{eq:2dflux}, we can compute $h$ once we have the
approximation for these normal fluxes.  In other words, we can solve the $1$D
Riemann problems for piecewise constant jumps locally in $x$ and $y$
directions, then sum these fluxes across the cell boundary $\xi_{i,j+\half}$ to
obtain the flux between hyperplanes.  One thereby decomposes the
multi-dimensional Riemann problem into a set of single-dimensional ones, to be
combined together by the inversion formula \eqref{eq:inverse}.

Unlike in the constant coefficient case, the flux function $h(q)$ must be
updated at every time step, as is usually done for finite volume methods,
although one may apply the nonlinear LTS method on the transformed problem
regardless. This would be based on the 1D analogues studied in
\cite{largetimestep1,largetimestep2, largetimestep3}.  The fully nonlinear
splitting will not be implemented here, but will be investigated in a
future work.

\section{Discrete Radon transform (DRT)}\label{sec:drt}

There are many different discretizations of the Radon transform and its inverse
\cite{matus, beylkin,brandt,kelley,framework-drt}, arguably the most well-known
being the filtered backprojection (FBP) algorithm \cite{ctbook}.  However, its
reliance on Fourier transforms and spherical harmonics lead to some filtering
of high-frequency content, causing Gibbs phenomenon near the sharp edges in the
solution.  This is not suitable for use in hyperbolic PDEs, which are known to
develop shock discontinuities.

Instead, we consider the use of a completely discrete analogue, namely the
approximate discrete Radon transform (ADRT), which we refer to simply as the
discrete Radon transform (DRT), introduced in \cite{GD96,brady}.  Rather than
interpolating pixel values onto straight lines passing near it, DRT sums one
entry for each row or column, along so-called \emph{digital lines} or
\emph{d-lines}. The \emph{d-lines} are defined recursively, allowing for a fast
computation in $\cO(N^2\log N)$ for an image of size $N \times N$.  The
back-projection is given by reversing the recursion, and is also fast with the
same computational cost of $\cO(N^2\log N)$.  The precise definitions are given
below.

\subsection{Recursive definition of DRT}

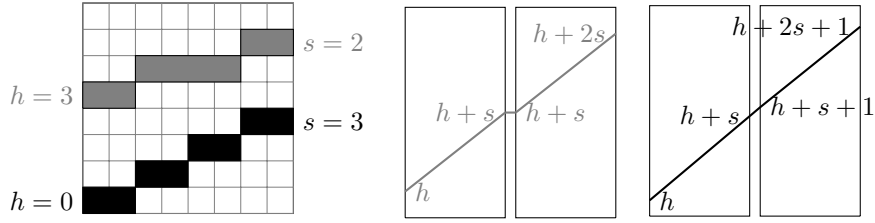
\begin{figure}
\centering
\begin{tabular}{ccc}
\begin{tikzpicture}[scale=0.35]
\draw (0,0) -- (8,0) -- (8,8) -- (0,8) -- (0,0);
\draw[step=1,gray,very thin] (0,0) grid (8,8);
\draw[fill=black] (0,0) -- (2,0) -- (2,1) -- (0,1) -- (0,0);
\draw[fill=black] (2,1) -- (4,1) -- (4,2) -- (2,2) -- (2,1);
\draw[fill=black] (4,2) -- (6,2) -- (6,3) -- (4,3) -- (4,2);
\draw[fill=black] (6,3) -- (8,3) -- (8,4) -- (6,4) -- (6,3);
\draw[black] (0,0.5) node[anchor=east] {$h=0$};
\draw[black] (8,3.5) node[anchor=west] {$s=3$};
\draw[fill=gray] (0,4) -- (2,4) -- (2,5) -- (0,5) -- (0,4);
\draw[fill=gray] (2,5) -- (6,5) -- (6,6) -- (2,6) -- (2,5);
\draw[fill=gray] (6,6) -- (8,6) -- (8,7) -- (6,7) -- (6,6);
\draw[gray] (0,4.5) node[anchor=east] {$h=3$};
\draw[gray] (8,6.5) node[anchor=west] {$s=2$};
\end{tikzpicture}
&
\begin{tikzpicture}[scale=0.7]
\draw (-2,0) -- (-0.1,0) -- (-0.1,4) -- (-2,4) -- (-2,0); 
\draw (2,0) -- (0.1,0) -- (0.1,4) -- (2,4) -- (2,0); 
\draw[gray,thick] (-2,0.5) node[anchor=west] {$h$} -- (-0.1,2.0) node[anchor=east] {$h+s$} ;
\draw[gray,thick] (-0.1,2.0) -- (0.1,2.0);
\draw[gray,thick] (2,3.5) node[anchor=east] {$h+2s$} -- (0.1,2.0) node[anchor=west] {$h+s$};
\end{tikzpicture} 
&
\begin{tikzpicture}[scale=0.7]
\draw (-2,0) -- (-0.1,0) -- (-0.1,4) -- (-2,4) -- (-2,0); 
\draw (2,0) -- (0.1,0) -- (0.1,4) -- (2,4) -- (2,0); 
\draw[black,thick] (-2,0.3) node[anchor=west] {$h$} -- (-0.1,1.9) node[anchor=east] {$h+s$} ;
\draw[black,thick] (2,3.6) node[anchor=east] {$h+2s+1$} -- (0.1,2.07) node[anchor=west] {$h+s+1$} -- (-0.1,1.9);
\end{tikzpicture} 
\end{tabular}
\caption{Examples of digital lines (d-lines) determined by two parameters 
    $h$ and $s$ (left) and the diagram of the recursion relation
\eqref{eq:dline} (right). In both figures the case when $s$ is even is in 
gray, and the case $s$ is odd is in black.} 
\label{fig:dlines}
\end{figure}

The \emph{d-lines} of length $N$ are denoted by $D_N(h,s)$ 
with two parameters $h$ and $s$ (see Figure \ref{fig:dlines}.)
$h$ denotes the \emph{height} ($x$-intercept) and 
$s$ the \emph{slope} ($x$-displacement), and the pair 
corresponds to $s$ and $\boldsymbol{\omega}$ for the continuous transform 
\eqref{eq:rt}, respectively.
Although the same notation $s$ is used here again after having been used in the 
continuous setting \eqref{eq:radon} we will keep the notation in order to follow the 
intuitive notation of \cite{press}, and mark the continuous variable with 
a subscript $s_c$ whenever the two are used simultaneously.
The definition uses the recursion in which d-lines of length $2n$ are split
into left and right d-lines of half its length,
\begin{equation}
    \begin{aligned}
    D_{2n}(h,2s) &= D^L_{n}(h,s) \cup D^R_{n}(h+s,s), \\
    D_{2n}(h,2s+1) &= D^L_{n}(h,s) \cup D^R_{n}(h+s+1,s).
    \end{aligned}
    \label{eq:dline}
\end{equation}
The recursion \eqref{eq:dline} defines only a quarter of the possible d-lines,
as the slope $s$ will range from 0 to $N$, corresponding to angles $0$ to
$\pi/4$ starting from the $x$-axis in the counter-clockwise direction.  This is
referred to by saying that the d-lines cover one quadrant, for the full
transform one needs to cover the angles from $0$ to $\pi$.  The other d-lines
can be computed by transposing and flipping the indices $h$ and $s$.  We will
denote the d-lines and DRT corresponding to the angular intervals
$[0,\pi/4]$,$[\pi/4,\pi/2]$,$[\pi/2,3\pi/4]$, and $[3\pi/4,\pi]$ by
$a$,$b$,$c$, and $d$.

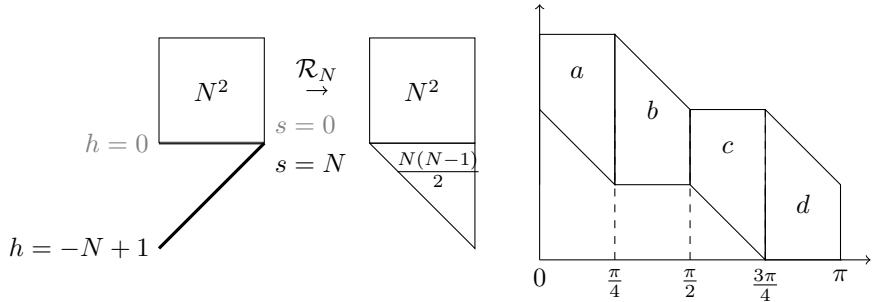
\begin{figure}
\centering
\begin{tabular}{cc}
\begin{minipage}{0.5\textwidth}
\begin{tikzpicture}[scale=0.7]
\draw (1,1) node {$N^2$};
\draw[gray,very thick] (-0,0) node[anchor=east]  {$h=0$};
\draw[gray,very thick] (-0,0) -- (2,0) node[anchor=south west] {$s=0$};
\draw[black] (-0,-2) node[anchor=east]  {$h=-N+1$};
\draw[black,very thick] (-0,-2) -- (2,0) node[anchor=north west] {$s=N$};
\draw (0,0) -- (2,0) -- (2,2) -- (0,2) -- (0,0);
\draw[->] (2.8,1) -- (3.0,1) node[anchor=south]  {$\mathcal{R}_N$} -- (3.2,1);
\draw (4,0) -- (6,0) -- (6,2) -- (4,2) -- (4,0);
\draw (5,1) node {$N^2$};
\draw (4,0) -- (6,0) -- (6,-2) -- (4,0);
\draw (5.3,-0.5) node {$\frac{N(N-1)}{2}$};
\end{tikzpicture} 
\end{minipage}
&
\begin{minipage}{0.4\textwidth}
\begin{tikzpicture}
\draw[dashed] (1,-2) node[anchor=north] {$\frac{\pi}{4}$} -- (1,0);
\draw[dashed] (2,-2) node[anchor=north] {$\frac{\pi}{2}$} -- (2,0);
\draw[dashed] (3,-2) node[anchor=north] {$\frac{3\pi}{4}$} -- (3,0);
\draw[dashed] (4,-2) node[anchor=north] {$\pi$} -- (4,-2);
\draw (0,0) --  (1,-1) -- (1, 1) -- (0, 1) -- (0,-1);
\draw (1,-1) -- (2,-1) -- (2, 0) -- (1, 1) -- (1,-1);
\draw (2,0) --  (2,-1) -- (3,-2) -- (3, 0) -- (2, 0);
\draw (3,0) --  (4,-1) -- (4,-2) -- (3,-2) -- (3,-1);
\draw (0,-2) node[anchor=north] {$0$};
\node at (0.5,0.5) {$a$};
\node at (1.5,0.0) {$b$};
\node at (2.5,-0.5) {$c$};
\node at (3.5,-1.25) {$d$};
\draw[->] (0,-2) -- (4.4,-2) node[anchor=north] {} ;
\draw[->] (0,-2) -- (0,1.4) node[anchor=east] {};
\end{tikzpicture}
\end{minipage}
\end{tabular}
\caption{The range of a quadrant of a discrete Radon transform (left)
and a diagram showing how the boundary of the quadrants $\{a,b,c,d\}$ 
can be identified (right). Here $\theta = \arctan(s/(N-1))$.}
\label{fig:drt_range}
\end{figure}

\begin{figure}
    \centering
    \begin{tabular}{rc}
        \rotatebox{90}{$\qquad \qquad \qquad t = 0$} &
        \includegraphics[width=0.6\textwidth]{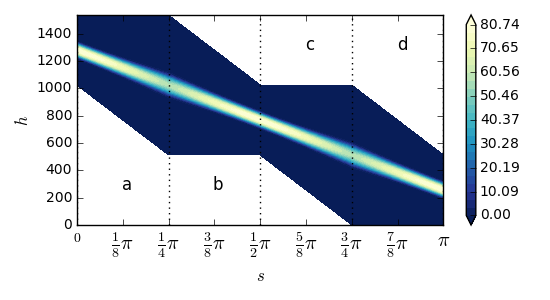}\\
        \rotatebox{90}{$\qquad \qquad \qquad t = 1$} &
        \includegraphics[width=0.6\textwidth]{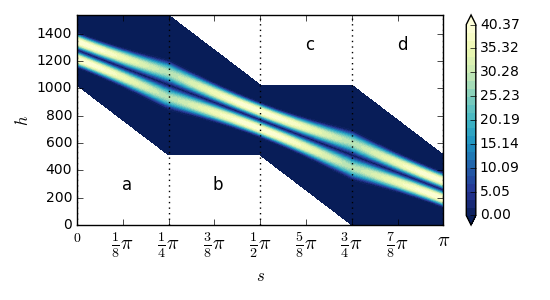}\\
        \rotatebox{90}{$\qquad \qquad \qquad t = 3$} &
        \includegraphics[width=0.6\textwidth]{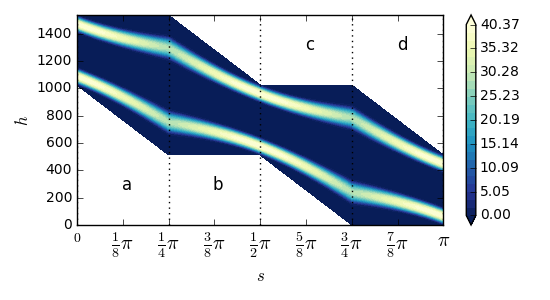}\\
    \end{tabular}
    \caption{The discrete Radon transform (DRT) of the solution to the 
    acoustic equation \eqref{eq:acoustic2d}, for times $t=0,1$ and $3$.
    The parameters $h$ and $s$ which appear on the axes 
    designate \emph{d-lines} (see Figure \ref{fig:dlines})
    and indices $\{a,b,c,d\}$ denote quadrants (see Figure \ref{fig:drt_range}).
    Details appear in Section \ref{sec:drt}. For a comparison with the 
    continuous Radon transform, see the right column of 
    Figure \ref{fig:acoustic}.}
    \label{fig:acoustic_drt}
\end{figure}

\begin{figure}
    \centering
    \begin{tabular}{rc}
   \rotatebox{90}{$\qquad \qquad \qquad t = 0$} &
   \includegraphics[width=0.6\textwidth]{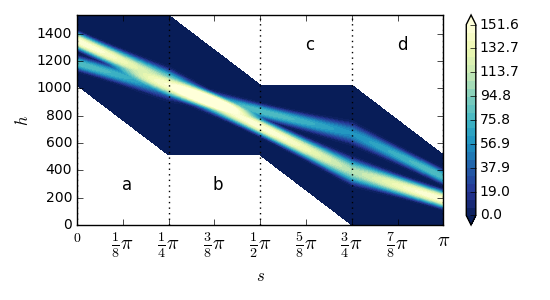}\\
   \rotatebox{90}{$\qquad \qquad \qquad t = 1.0$} &
   \includegraphics[width=0.6\textwidth]{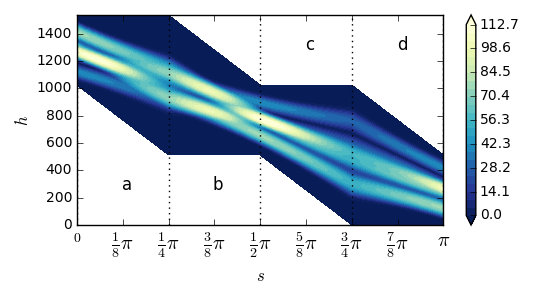}\\
   \rotatebox{90}{$\qquad \qquad \qquad t = 1.5$} &
   \includegraphics[width=0.6\textwidth]{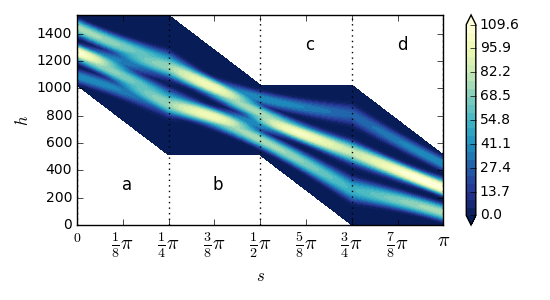}\\
    \end{tabular}
    \caption{The discrete Radon transform (DRT) of the solution to the acoustic
equation \eqref{eq:acoustic2d} with initial conditions \eqref{eq:2humps}, for
times $t=0,1$ and $1.5$.  The parameters $h$ and $s$ which appear on the
axes designate \emph{d-lines} (see Figure \ref{fig:dlines}) and indices
$\{a,b,c,d\}$ denote quadrants (see Figure \ref{fig:drt_range}). For a
comparison with the continuous Radon transform, see the right column of Figure
\ref{fig:acoustic_2humps}.}
    \label{fig:acoustic_2humps_drt}
\end{figure}

The \emph{DRT of an array $A \in \RR^2$ for the quadrants $a,b,c$ and $d$} are
given by the summation of entries of $A$ over the d-lines,
\begin{equation}
    \begin{aligned}
    \left(\cR_N^a A\right)_{h,s} &= \sum_{(i,j) \in D_N(h,s)} A_{i,j}, \\
    \left(\cR_N^b A\right)_{h,s} &= \sum_{(i,j) \in D_N(h,s)} A_{j,i}, \\
    \left(\cR_N^c A\right)_{h,s} &= \sum_{(i,j) \in D_N(h,s)} A_{j,N-i+1}, \\
    \left(\cR_N^d A\right)_{h,s} &= \sum_{(i,j) \in D_N(h,s)} A_{N-i+1,j}. \\
    \end{aligned}
    \label{eq:drtsum}
\end{equation}
The full DRT is simply the ordered tuple of all quadrants, and we write
\begin{equation}
\cR_NA = (\cR_N^a A, \cR_N^b A, \cR_N^c A, \cR_N^d A).
\end{equation}
See \cref{fig:drt_range} for a visual illustration. Due to the recursive form of
\eqref{eq:dline}, the transform can be computed in $\cO(N^2 \log N)$. The
parameters $h$ and $s$ belong to the range $[-s+1,N]$ and $[0,N]$ so $R^a_N A
\in \RR^{\left(\frac{3}{2}N + \frac{1}{2} \times N\right)}$.  Therefore,
$\cR_N: \RR^{N \times N} \to \RR^{(6N + 2) \times N}$. For example, the DRT of
the solutions displayed in \cref{fig:acoustic,fig:acoustic_2humps} are plotted
in \cref{fig:acoustic_drt,fig:acoustic_2humps_drt}, respectively.

There is a simple relationship between the DRT and the continuous Radon
transform. First let us say that $s_c \in [-1,1]$ (perhaps through proper
scaling) and parameterize $\omega$ by $\omega = (\cos \theta, \sin \theta)$
where $\theta \in [0,\pi]$ .  The relation to continuous variables
$(s_c,\omega)$ is given by
\begin{equation}
    s_c = \cos \theta\left(\frac{2h}{N} - 1 + \frac{s}{N-1}\right),
        \quad \quad
        \theta = \arctan\left( \frac{s}{N-1} \right).
\end{equation}
Then the explicit relation between the DRT $\cR_N$ and the continuous transform
$\cR$ are given after the density of the lines are also transformed depending
on the angle by $\cos \theta$, 
\begin{equation} 
    \begin{aligned}
    \cR_N^a f (h, s) &= \cos \theta \, \cR(s_c,\theta),\\
    \cR_N^b f (h, s) &= \cos \theta \, \cR(s_c,\pi - \theta),\\
    \cR_N^c f (h, s) &= \cos \theta \, \cR(s_c,3\pi/2 - \theta),\\
    \cR_N^d f (h, s) &= \cos \theta \, \cR(s_c,3\pi/2 + \theta).
    \end{aligned}
    \label{eq:drt2rt}
\end{equation}
We note that the DRT approximates the continuous transform with first order
accuracy with respect to the grid cell width $1/N$ \cite{brady}.

The back-projection is the dual of this transform with respect to the usual dot
product in $\RR^{N^2}$.  We will denote the back-projection by $B_N$ or
$\cR_N^T$. If one explicitly forms the matrix for the linear transforms $\cR_N$
and $\cR_N^T$ they are indeed transposes of each other.  $\cR_N^T$ is the
discrete analogue of $\cR^\#$ in \eqref{eq:backprojection}, a summation of all
values assigned to d-lines passing through a point.

$\cR^T_N$ is computed by reversing the sweep \eqref{eq:dline} above and
computing a sequence of back-projections of decreasing size.  Given a matrix
$\hat{A} \in \RR^{\left(\frac{3}{2}N + \frac{1}{2}\right) \times N}$, the
reverse sweep for one quadrant is given by \begin{equation} \begin{aligned}
B^L_n(h,s) &=  B_{2n}(h,2s) + B_{2n} (h,2s+1),\\ B^R_n(h+s,s) &=  B_{2n}(h,2s)
+ B_{2n}(h-1,2s), \end{aligned} \end{equation} where the initial array
$B_N(h,s) = \hat{A}_{h,s}$, and $n$ is set to $N/2$.  This summation is
repeated for the two half-images $B_n^L$ and $B_n^R$ on the LHS, until $n$
reaches $1$. Again, this summation is only for one quadrant, and we denote the
end result as $\left(B_N^a \hat{A} \right)_{i,j}$.

The full back-projection is given by
\begin{equation}
    \left(\cB_N \hat{A} \right) (i,j) = 
    \frac{1}{4N^2} \left( B_N^a \hat{A} + B_N^b \hat{A} 
                        + B_N^c \hat{A} + B_N^d \hat{A} \right).
\end{equation}

Its inversion algorithm using a full multi-grid method was demonstrated 
in \cite{press} along with convergence analysis. 
In this paper, we use the conjugate gradient (CG) algorithm 
as will be discussed below in Section \ref{sec:idrt}.

We end this section with the remark that the recursion \eqref{eq:dline}
need not be in two and can be in any prime number, much like the 
fast Fourier transform \cite{cooley-tukey}.

\subsection{Inversion of DRT with Conjugate Gradient Method}
\label{sec:idrt}

In order to use the dimensional splitting method to solve PDEs, a method for
computing the inverse of a DRT \eqref{eq:inverse} is needed.  An inversion
algorithm using a full multi-grid method appeared in \cite{press}. Here we
explore the application of the conjugate gradient method \cite{greenbaum} to
the least-squares problem
\begin{equation}
    \cR_N^T \cR_N X = \cR_N^T B.
    \label{eq:normaleqn}
\end{equation}
The matrices for the transforms $\cR^T_N$ and $\cR_N$ are never explicitly
formed, as we can use the fast algorithm.  The computational cost of a DRT
inversion is conjectured to be $\cO(N^{5/2} \log N)$ for an $N \times N$
image \cite{rimThesis}.  Note that this is slightly more costly than
$\cO(N^2 (\log N)^3)$ that was conjectured for the full multi-grid method
\cite{press}.  A more careful study of this inversion is of interest on its own
right, and will appear elsewhere.

The inversion of the Radon transform, be it continuous or discrete, is
mildly ill-posed \cite{ctbook}. Numerically speaking, this means the matrix
$\cR^T_N \cR_N$ operator to be inverted will be ill-conditioned. One approach
commonly used to deal with this issue is to use regularizations, for example in
medical imaging applications. However, there is an important distinction to the
tomography setting, namely that it is feasible to make additional measurements.
In our setting, making more measurements from the original image $X$ in
\eqref{eq:normaleqn} would only incur additional \emph{computational} effort,
whereas in medical imaging it would require more \emph{physical} measurements.

For example, in this work we perturb the range of the Radon transform and
there is the possibility the perturbed function on the space of hyperplanes no
longer lies in the range of the transform. The inversion \eqref{eq:normaleqn}
is exact only when $B$ lies in the range of $\cR_N$, and this assumption cannot
be satisfied in general once $B$ is evolved with respect to the dynamics of the
transformed variables, as in \eqref{eq:solveRu}.  Therefore, changes in the
transformed variables will cause $\hat{q}$ to depart from the range of $\cR_N$.
This becomes a source of error, incurring numerical artifacts in the computed
inverse. The DRT employed in this paper also does incur these artifacts.

A brute force way to avoid this error is to make more measurements,
i.e., \emph{oversample}. We simply prolong the original image $q$ before
manipulation, and restrict after the back-projection. This would correspond
to making additional measurements in the tomography setting. This enlarges the
range of the transform, and allows one to control the error. Therefore $\cR_N$
will be replaced by $\cR_{2pN} \cP_{2p}$ where $\cP_{p}$ is the $0$-th order
prologation (where the value of each cell in the original grid is assigned to a
$2p \times 2p$ cells in the enlarged grid) and $\cR^T_N$ by $\cS_{2p}
\cR^T_{2pN}$ where $\cS$ is the restriction operator. The oversampling
strategy does not affect the overall complexity.

\subsection{DRT in dimension three} \label{sec:drt3}

Just as the continuous Radon transform was defined in \eqref{eq:radon} for
arbitrary number of dimensions $n$, the DRT can also be generalized to higher
dimensions \cite{nDRT}.  Here we treat the 3D case as an example.  The
recursive definitions \eqref{eq:dline} for the \emph{d-planes} parametrized by
three parameters $(h,s_1,s_2)$ can be derived for each hexadecant in a
straightforward manner, as follows
\begin{equation}
    \begin{aligned}
    D_{2n}(&h,2s_1,   2s_2)=            
                                D^{LL}_n(h        ,s_1,  s_2) 
                           \cup D^{RL}_n(h+s_1    ,s_1,  s_2) \\
             & \quad \quad \cup D^{LR}_n(h+s_2    ,s_1,  s_2)
                           \cup D^{RR}_n(h+s_1+s_2,s_1,  s_2),\\
    D_{2n}(&h,2s_1+1, 2s_2)= 
                                D^{LL}_n(h          ,s_1,  s_2) 
                           \cup D^{RL}_n(h+s_1+1    ,s_1,  s_2) \\
             & \quad \quad \cup D^{LR}_n(h+s_2      ,s_1,  s_2)
                           \cup D^{RR}_n(h+s_1+s_2+1,s_1,  s_2),\\
    D_{2n}(&h,2s_1,   2s_2+1)= 
                                D^{LL}_n(h          ,s_1,  s_2) 
                           \cup D^{RL}_n(h+s_1      ,s_1,  s_2) \\
             & \quad \quad \cup D^{LR}_n(h+s_2+1    ,s_1,  s_2)
                           \cup D^{RR}_n(h+s_1+s_2+1,s_1,  s_2),\\
    D_{2n}(&h,2s_1+1, 2s_2+1)=
                                D^{LL}_n(h          ,s_1,  s_2) 
                           \cup D^{RL}_n(h+s_1+1    ,s_1,  s_2) \\
             & \quad \quad \cup D^{LR}_n(h+s_2+1    ,s_1,  s_2)
                           \cup D^{RR}_n(h+s_1+s_2+2,s_1,  s_2).\\
    \end{aligned}
\end{equation}
The DRT over one hexadecant (a quarter of a quadrant) is defined as the sum over 
these d-planes as in \eqref{eq:drtsum}, and now the full transform in 3D is given by 
applying these to each of the hexadecant 
\begin{equation}
    \cH = 
    \left\{
    \begin{matrix*}[r]
          aa,& ab,& ac,& ad, \\
          ba,& bb,& bc,& bd, \\
          ca,& cb,& cc,& cd, \\
          da,& db,& dc,& {dd\,}
    \end{matrix*}
    \right\}.
\end{equation}
Hence, via transposing and flipping the indices as necessary, the full DRT is
the ordered tuple
\begin{equation}
        \cR_NA = 
        \left(
        \begin{matrix*}[r]
          \cR_N^{aa} A,& \cR_N^{ab} A,& \cR_N^{ac} A,& \cR_N^{ad} A,\\
          \cR_N^{ba} A,& \cR_N^{bb} A,& \cR_N^{bc} A,& \cR_N^{bd} A,\\
          \cR_N^{ca} A,& \cR_N^{cb} A,& \cR_N^{cc} A,& \cR_N^{cd} A,\\
          \cR_N^{da} A,& \cR_N^{db} A,& \cR_N^{dc} A,& {\cR_N^{dd} A\,}
        \end{matrix*}
        \right)
\end{equation}

The corresponding back-projection operation for a hexadecant is given by
\begin{equation}
    \begin{aligned}
        B_n^{LL}&(h,s_1, s_2)=            
                                 B_{2n}(h, 2s_1  , 2s_2  ) 
                               + B_{2n}(h, 2s_1+1, 2s_2  ) \\
                   &           + B_{2n}(h, 2s_1  , 2s_2+1)
                               + B_{2n}(h, 2s_1+1, 2s_2+1),\\
       B_n^{RL}&(h+s_1,s_1, s_2)= 
                                 B_{2n}(h,  2s_1,  2s_2  ) 
                               + B_{2n}(h-1,2s_1+1,2s_2  ) \\
                   &           + B_{2n}(h,  2s_1,  2s_2+1)
                               + B_{2n}(h-1,2s_1+1,2s_2+1),\\
       B_n^{LR}&(h+s_2,s_1, s_2)= 
                                 B_{2n}(h  ,2s_1,  2s_2  ) 
                               + B_{2n}(h  ,2s_1+1,2s_2  ) \\
                   &           + B_{2n}(h-1,2s_1,  2s_2+1)
                               + B_{2n}(h-1,2s_1+1,2s_2+1),\\
       B_n^{RR}&(h+s_1+s_2,s_1, s_2)=
                                 B_{2n}(h  ,2s_1,  2s_2  ) 
                               + B_{2n}(h-1,2s_1+1,2s_2  ) \\
                   &           + B_{2n}(h-1,2s_1,  2s_2+1)
                               + B_{2n}(h-2,2s_1+1,2s_2+1).\\
    \end{aligned}
\end{equation}
The full back-projection is then the average of back-projections $B_N$ over all
hexadecants in $\cH$,
\begin{equation}
    \left(B_N \hat{A} \right) (i,j) =  
    \frac{1}{16 N^3} \sum_{k \in \cH} B_N^k \hat{A}.
\end{equation}
The computational cost of both operations would be $\cO(N^3 \log N)$.

\section{Applications of the dimensional splitting}\label{sec:radonsplituse}

The dimensional splitting described in Section \ref{sec:radonsplit} above is a
decomposition of hyperbolic solutions into evolution of planar waves. This
decomposition can be useful in diverse settings. Here we discuss two
applications: the absorbing boundary conditions and the displacement
interpolation.

\subsection{Absorbing boundary conditions}\label{sec:absorb}

It is well-known that imposing absorbing boundary conditions to emulate
infinite domains in multi-dimensional wave propagation is a challenging problem
\cite{abc,pml,stretched-coords,pml3d}.  On the other hand, the $1$D
extrapolation boundary condition is much more tractable \cite{fvmbook}. A
major advantage of this splitting method is that the 1D extrapolation boundary
conditions can be used on the transformed side at the computational boundary to
avoid any reflections.  This yields \emph{exactly} the desired absorbing
boundary conditions in the odd-dimensional case. Therefore, the dimensional
splitting introduced in the previous section can be used directly to impose
absorbing boundary conditions in 3D.  On the contrary, there is an error caused
by such an extrapolation in the even-dimensional case. This is due to the
Huygens' principle, evident in the presence of the Hilbert transform in the
inversion formula \eqref{eq:inverse}.  In this section, we discuss the type of
error caused by imposing such extrapolation boundary conditions via the Radon
transform in even dimensions.

In the true infinite domain, the non-zero values in transformed variables
beyond the computational boundary of $\cS^{n-1} \times \RR$ affect the solution
within the computational domain in the original variables $\RR^n$.  For
example, the vertical translation of horizontal strips in Figure
\ref{fig:acoustic} should continue beyond the finite computational boundary,
and by imposing a 1D extrapolation boundary condtion we would be neglecting
this infinite propagation.  To make this more precise, denote the computational
(finite) transformed domain by $\Omega = \{(\omega,s) \in S^{n-1} \times
(-b,b)\}$ for some $b>0$.  Let $\chi_\Omega$ be the characteristic function of
the finite domain and $\chi_{\RR^n \setminus \Omega} = 1- \chi_\Omega$.
For $n$ even, the exact solution $q$ can be written as,
\begin{align}
  q(T,x)  
  &= \frac{1}{c_n}\cR^\# H_s \frac{d^{n-1}}{ds^{n-1}} \, \hat{q}(T,\omega,s) \\
  &= \frac{1}{c_n}\cR^\# H_s \chi_\Omega \frac{d^{n-1}}{ds^{n-1}} \, 
   \hat{q}(T,\omega,s) 
   +  \frac{1}{c_n}\cR^\# H_s \chi_{\RR^n \setminus \Omega}
   \frac{d^{n-1}}{ds^{n-1}} \, \hat{q}(T,\omega,s).
   \label{eq:extrapsum}
\end{align}
Recall that $\cR^\#$ is the back-projection \eqref{eq:backprojection}.  The
first term in \eqref{eq:extrapsum} is the approximate solution one would obtain
if extrapolation boundary was set up at the boundaries $s = \pm b$.  Let us
call this approximate solution $q_h(x)$. Then the error is
\begin{align}
    q(T,x)  - q_h(T,x) 
         &=  \frac{1}{c_n}\cR^\# 
         \left(\mathrm{p.v}
         \int_{(-\infty,-b)] \cup [b, \infty)} 
         \frac{1}{z-s}
         \frac{\p^{n-1}}{\p z^{n-1}} \, 
     \hat{q}(T,\omega,z)\mathrm{d}z \right),
            \label{eq:absorberror0}
\end{align}
where $\mathrm{p.v}$ denotes the principal value integral.  Note that in
hyperbolic problems in free space, wave profiles will radiate outwards, that
is, the support of $\hat{q}$ will be transported towards $r = \pm \infty$. This
causes the RHS above to decay with time. Furthermore, the principal value
integral is a smooth function of $s$ as long as $\p^{n-1} \hat{q} / \p r^{n-1}$
is integrable. Since we will also apply $\cR^\#$, we expect the error to be
smoother than $q(x)$.

Let us revisit the acoustic equations example \eqref{eq:acoustic2d} from
Section \ref{sec:linearlts}, with initial conditions \eqref{eq:coshump}.  Since
$\hat{q}_0(x_1,x_2)$ is supported in  $\{(\omega,s) \in S^1
\times \RR: \abs{s} \le 1 \}$, we have a simple estimate for the case when $t$
is sufficiently large so that $(-1+ t,1+ t) \subset (b,\infty)$, 
\begin{equation}
\begin{aligned}
    \Norm{q(t,x) - q_h(t,x)}{1} &\le \frac{1}{2c_2} \cR^\#
    \left\lVert
        \int_{-\infty}^{-b}  
         \frac{r_1}{z-s} 
         \frac{d}{dz} \, 
           \hat{p}_0(z + t) 
          \, \mathrm{d}z \right. \\
          &\qquad \qquad \qquad \qquad \left.+ 
            \int_b^\infty
         \frac{r_2}{z-s} 
            \frac{d}{dz} \, 
             \hat{p}_0(z - t) 
            \, \mathrm{d}z \right\rVert_1 \\
            & \le 
            \frac{2\sqrt{6}\pi}{c_2 \abs{t - \Norm{x}{2} - 1}}
            \int_\RR \abs{\frac{d\hat{p}_0}{dz}} \, \mathrm{d}z,
\end{aligned} \label{eq:absorberror}
\end{equation}
where $\lVert \cdot \rVert_1$ is the $\ell^1$-norm for $\RR^3$.
Therefore we see that the effect of the extrapolation boundary decays
relatively slowly, at the rate of $\cO(1/t)$. 

\begin{figure}
\centering
\includegraphics[width=0.8\textwidth]{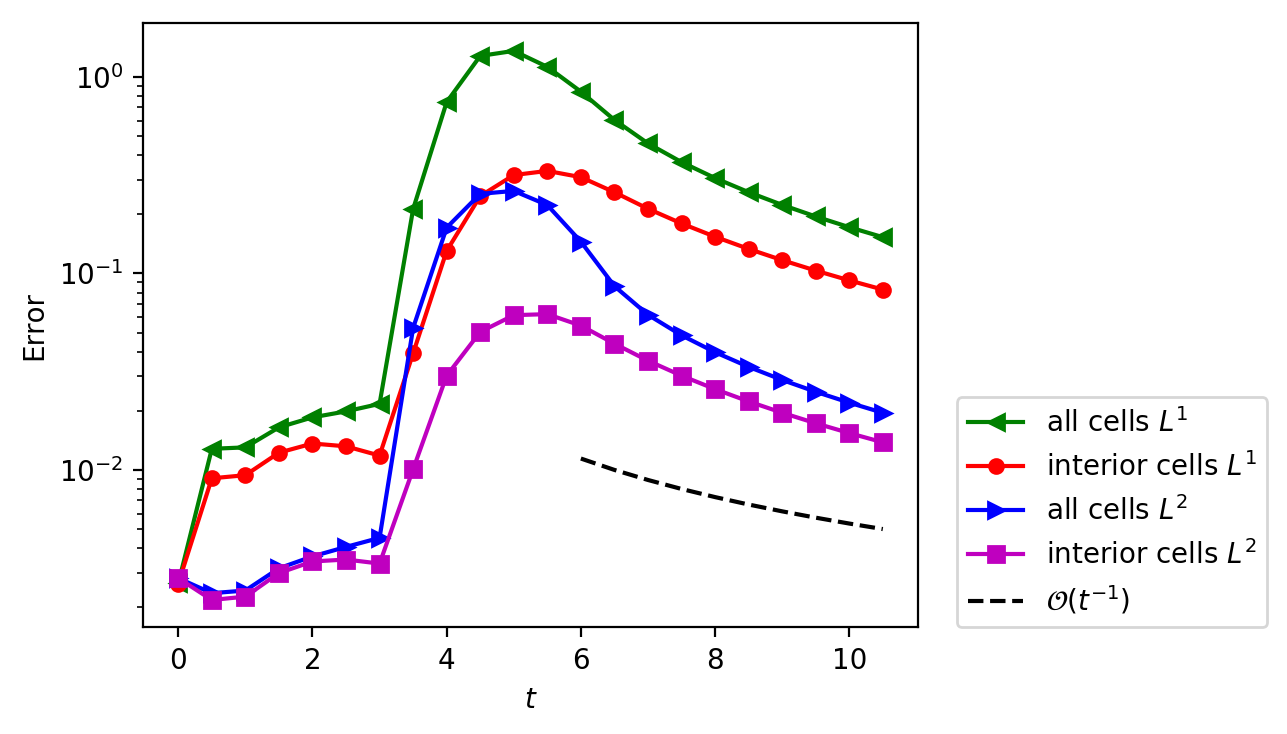}
\caption{$L^1$- and $L^2$-norm difference between the splitting
solution and the reference solution over time, for the acoustic equation
\eqref{eq:acoustic2d} with the initial condition \eqref{eq:coshump}. The
difference over all cells are shown together with the difference over interior
cells inside $[-3,3] \times [-3,3]$.  The slope of $\cO(1/t)$ is also shown for
comparison.}
\label{fig:acousticerror}
\end{figure}

The solution is now compared with a fully 2D reference finite volume solution
computed on a larger domain $[-8,8] \times [-8,8]$, using the wave propagation
algorithm \cite{rjl:wpalg} with Lax-Wendroff flux and Van Leer limiter
\cite{vanleer79,fvmbook}, implemented in \textsc{Clawpack}.  The reference
solution was computed on a $1024 \times 1024$ finer grid-cells of uniform size,
then corresponding cells have been summed and compared with coarser cells of
the DRT solution.  The $L^1$ and $L^2$-norms of the difference over
time is displayed below in Figure \ref{fig:acousticerror}. 

\begin{figure}
    \centering
    \begin{tabular}{m{0.25cm} m{0.35\textwidth} m{0.525\textwidth}}
        \rotatebox{90}{$t = 3$} &
        \includegraphics[width=0.35\textwidth]{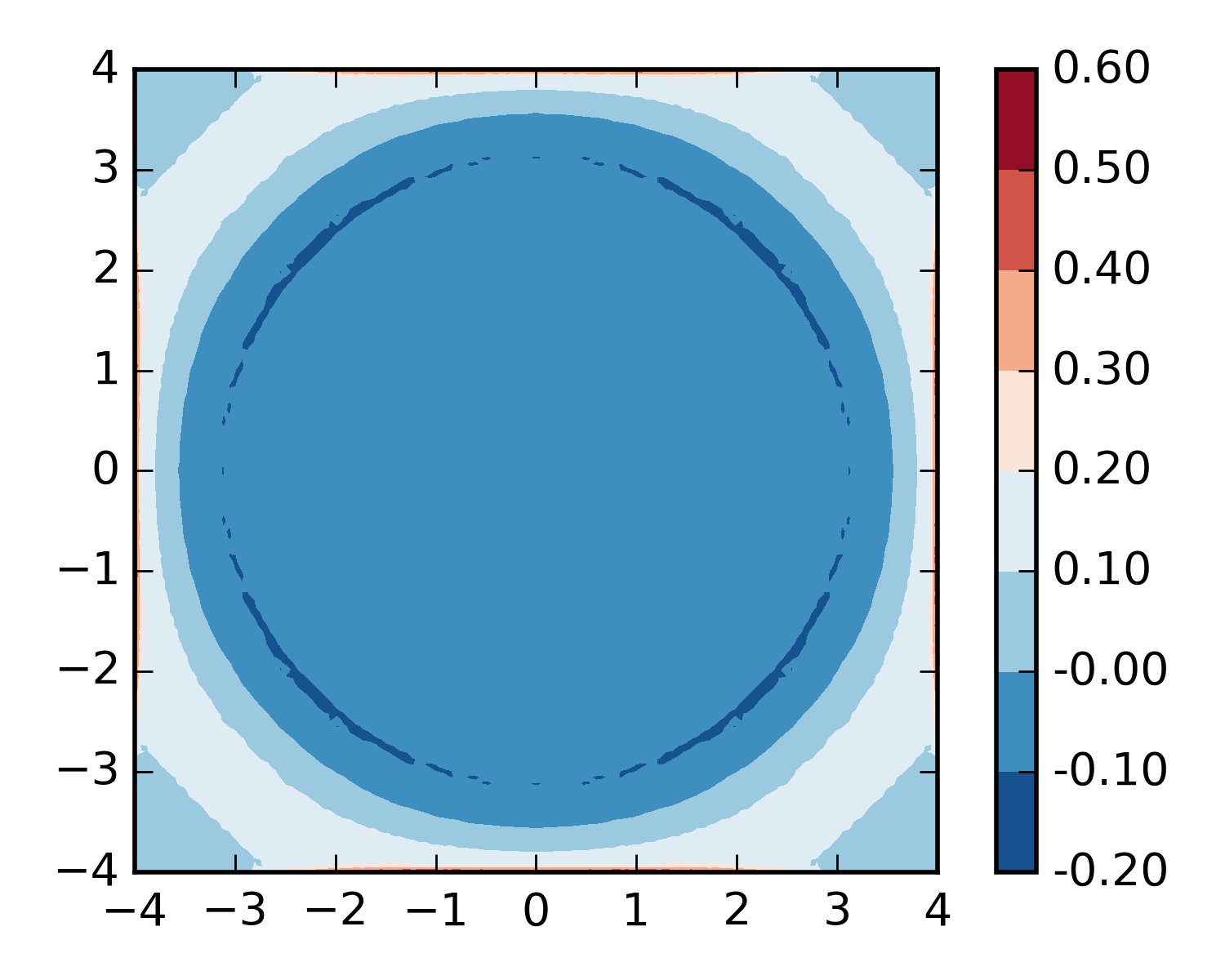}&
            \vspace{1.25em}
        \includegraphics[width=0.525\textwidth]{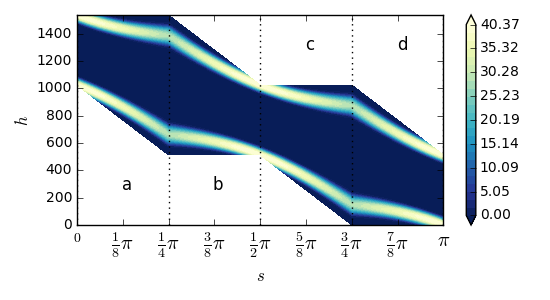}\\
    \end{tabular}
    \caption{The splitting solution to the acoustic equation \eqref{eq:acoustic2d} 
        in the square domain $[-4,4] \times [-4,4]$.
        The pressure $p$  and its DRT $\cR_N p$ at time $t=3$.}
        \label{fig:acousticbdry}
\end{figure}

The error is in the order of $10^{-3}$ up to time $t=2.5$, before the profile
starts approaching the boundary. The error from the truncation
\eqref{eq:absorberror} begins to appear around time $t=2.5$ and peaks around
time $t=5$, then decays to zero with time. The solution at time $t=3$, as it
has begun to interact with the boundary, is shown in Figure
\ref{fig:acousticbdry}.  Note that there are no reflections from this boundary
condition.  On the other hand, a thin layer appears at the computational
boundary.  The layer is clearly non-physical, but is localized and has limited
affect on the solution further in the interior. The DRT of the solution is also
shown to its right, and we can see that the two pulses from the d'Alembert
solution are hitting the 1D extrapolation boundaries (the top and bottom
boundaries of the polygonal region in dark blue). The pulses first arrive at
the DRT boundary at the angles $0, \pi/2$ and $\pi$, and this agrees with the
solution plot to the left.

\begin{figure}
    \centering
    \begin{tabular}{cc}
    \includegraphics[width=0.45\textwidth]{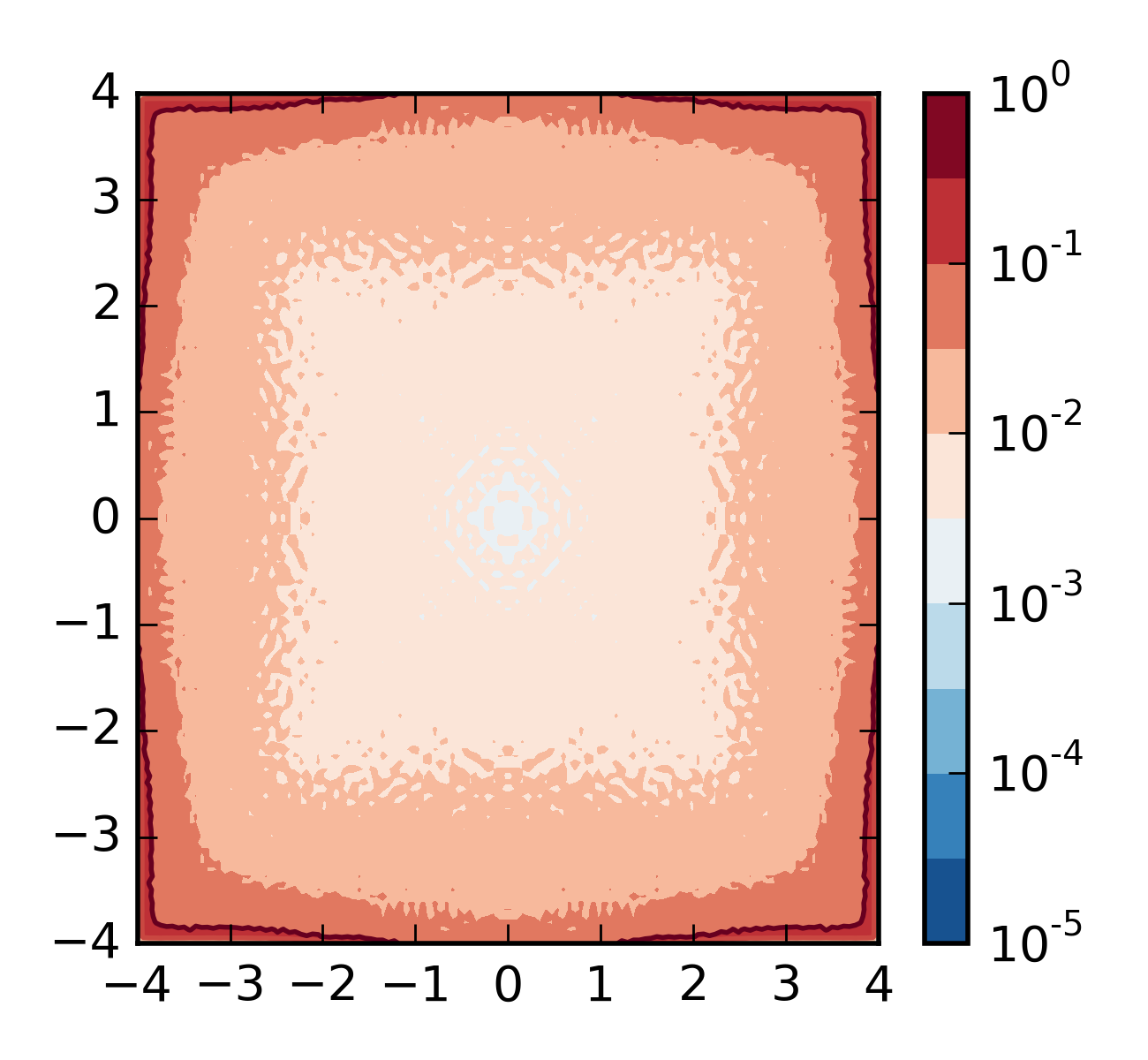}&
    \includegraphics[width=0.45\textwidth]{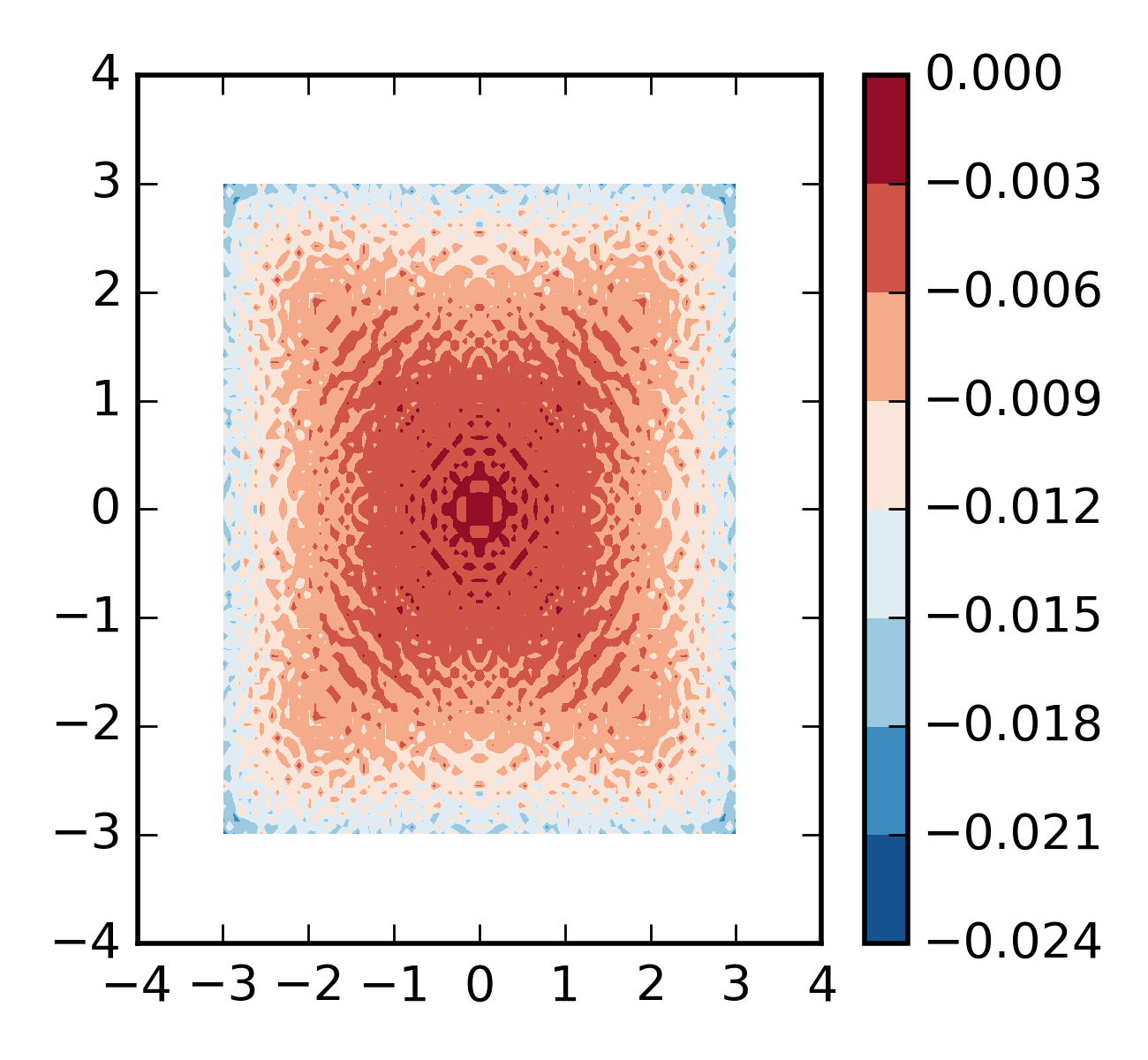}
    \end{tabular}
    \caption{Difference to the reference solution at time $t=5$.
             The difference for all cells in log-scale (left) and the difference
             for the interior cells in $[-3,3] \times [-3,3]$ (right).}
             \label{fig:interiorerror}
\end{figure}

Comparing this solution to the reference solution, one discovers that the bulk
of this error is concentrated near the thin layer which appears near the
boundary. In Figure \ref{fig:interiorerror}, we have plotted the difference
between our solution and the reference solution on the computational domain
$[-4,4] \times [-4,4]$ to the left.  When we restrict the contour plot to the
interior portion of the domain $[-3,3] \times [-3,3]$ as is shown to the right,
we see that the error is significantly smaller as we move away from the
boundary. The estimate \eqref{eq:absorberror} helps us understand
this behavior. Note the decay of the principal value integral in
\eqref{eq:absorberror0}: the further away the interior point is from the
boundary, the smaller is the effect of the trunctation by $\chi_\Omega$.
If we denote this distance by $d$, then the decay will be $\cO(1/(d+t))$.
As $t$ increases and the waves leave the domain, the error also decays, at 
the rate $\cO(1/t)$ estimated by \eqref{eq:absorberror}.

The observations above suggest several potential approaches in further
improvements to this approach. For example, one can exploit the fact that the
Hilbert transform in \eqref{eq:absorberror0} commutes with translation, to
emulate the effect of the infinite domain. One may also exploit the decay with
respect to the distance to the boundary by placing finer mesh along the
boundary in an adaptive fashion.

\subsection{Displacement interpolation}\label{sec:disp}

In projection-based model reduction, the solution to a parametrized PDE is
projected into a low-dimensional subspace, yielding a fast solver with
significantly lower computational cost without compromising accuracy.  To
discover this low-dimensional subspace, the popular approach is to use proper
orthogonal decomposition (POD) \cite{pod01}.  For hyperbolic PDEs, however, the
solutions do not lie in a low-rank linear subspace, even for the simplest
problems \cite{marsden1, amsallem, Carlberg15,Schulze15,reversal}. For
instance, the d'Alembert solution \eqref{eq:dalembert} is a linearly
independent function of $s$ for each $t > 0$.  It is easy to see that a linear
projection of this solution to a low-dimensional basis would not yield a good
approximation of the solution.  Naturally, methods to remove translational
symmetry \cite{marsden1,Schulze15,reversal} are being actively explored.  This
is also intimately related to the concept of \emph{displacement interpolation},
a term we borrow from the optimal transport literature
\cite{villani2008optimal}, in which one aims to minimize the Wasserstein
distance, although we will not make the connection more explicit here. 

\begin{figure}
    \centering
    \begin{tabular}{cc}
      \includegraphics[width=0.35\textwidth]{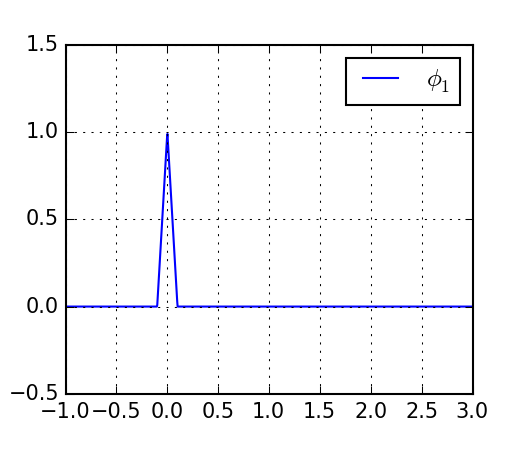}
        &
      \includegraphics[width=0.35\textwidth]{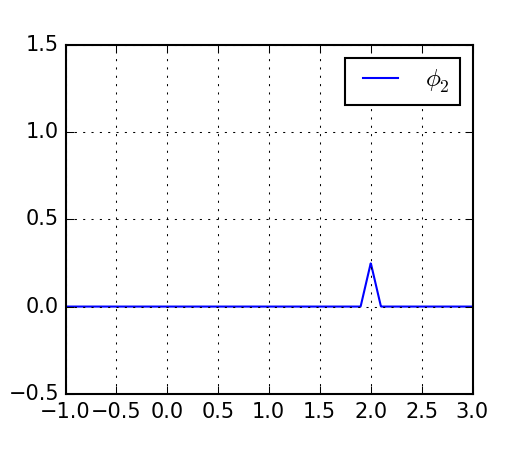}\\
      \includegraphics[width=0.35\textwidth]{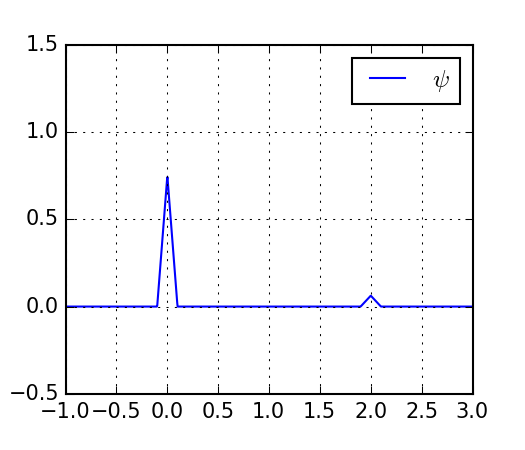}
        &
      \includegraphics[width=0.35\textwidth]{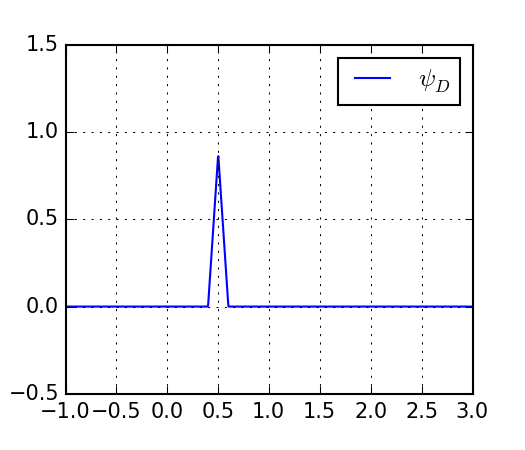}\\
    \end{tabular}
    \caption{Two hat functions $\phi_1$ and $\phi_2$ (top row) and 
    linear interpolation $\psi$ and 
    displacement interpolation $\psi_D$ between the two functions
    with respective weights $0.75$ and $0.25$ (bottom row).}
    \label{fig:disphat1d}
\end{figure}

Let us first illustrate how displacement interpolation arises naturally, with a
simple 1D example. Suppose $\phi_0$ is a hat function, given by
\begin{equation}
    \phi_0(x) = 
    \begin{cases} 
        \frac{x}{h} + 1 & \text{ if } -h < x < 0, \\ 
        -\frac{x}{h} + 1 & \text{ if } 0 \le x < h, \\
                              0 & \text{ otherwise,} \\
    \end{cases}
\end{equation}
    for some $h > 0$. Let $\phi_1$ and $\phi_2$ be translation and scaling of
    $\phi_0$,
\begin{equation}
    \phi_1(x) = \phi_0(x) 
    \quad \text{ and } \quad
    \phi_2(x) = \frac{1}{4}\phi_0(x - 2).
\end{equation}
For $h = 0.1$, the two functions are shown in the first row of Figure
\ref{fig:disphat1d}.  The linear interpolation $\psi$ of the two functions with
weights $(1-\tau)$ and $\tau$ is given by
\begin{equation}
\begin{aligned}
    \psi(x) &= (1 - \tau) \phi_1(x) + \tau \phi_2(x) \\
            &= (1 - \tau) \phi_0(x) + \frac{\tau}{4} \phi_0(x-2) 
\end{aligned}
\end{equation}
whereas a displacement interpolation between the two functions under a simple
transport map (1D translation) would be given by
\begin{equation}
\begin{aligned}
    \psi_D(x) &= (1 - \tau) \phi_1(x-2\tau) + \frac{\tau}{4}\phi_2(x + 2(1-\tau))  \\
              &= \left(1 - \frac{3}{4} \tau \right)\phi_0(x - 2\tau) 
\end{aligned}
\label{eq:1ddi}
\end{equation}
The two interpolants for $\tau = 0.25$ are plotted in the bottom row in Figure
\ref{fig:disphat1d} .  Since $\phi_1$ and $\phi_2$ are both translates of a
scalar multiple of $\phi_0$, the displacement interpolation reveals the
low-rank nature of the two functions, whereas the linear interpolant remains
rank two for $\tau \in (0,1)$.

In practice, one must be able to deduce that $\phi_1$ and $\phi_2$ above lie in
the translates of $\text{span}\{\phi_0\}$ without a priori knowledge.  To
achieve this, one may apply the \emph{template-fitting} procedure
\cite{marsden1} which solves the minimization problem
\begin{equation}
    \tau_* = \argmin_{\tau \in \RR} \Norm{\phi_2(x) - \cK(\tau)[\phi_1(x)]}{2},
    \label{eq:tf1d}
\end{equation}
where $\cK$ is the translation operator, $\cK(\tau) [\phi_1(x)] = \phi_1(x-
\tau)$, then perform a singular value decomposition (SVD) on $\{\phi_2,
\cK(\tau_*)\phi_1\}$ \cite{marsden1, Schulze15}.  However, this simple
formulation does not take into account multiple traveling speeds or heavily
deforming profiles, which limits its applicability. \emph{Transport reversal}
was introduced in \cite{reversal} to overcome these limitations.  The algorithm
is a greedy iteration over a generalized form of \eqref{eq:tf1d}, which
decomposes the 1D function $\phi_2(x)$ into multiple traveling structures. To
be more precise, given two functions $\phi_1$ and $\phi_2$ as in
\eqref{eq:tf1d}, transport reversal yields the decomposition
\begin{equation}
    \psi_D(x;\tau) = 
    \sum_{k=1}^K \eta_k(\tau) \cK(\nu_k\tau)[\rho_k(x,\tau)\varphi_1(x)].
    \label{eq:reversal1d}
\end{equation}
where $\eta_k$ is a scaling function and $\rho_k$ a cut-off function.  For more
detailed treatment of this decomposition in 1D, we refer the reader to
\cite{reversal}. Let us suppose we have computed this decomposition.  The
displacement interpolant $\psi_D$ resulting from this decomposition is set to
satisfy,
\begin{equation}
    \psi_D(x,0) = \phi_1(x) 
    \quad \text{ and } \quad
    \psi_D(x,1) = \phi_2(x).
\end{equation}

\begin{figure}
    \centering
    \begin{tabular}{cc}
        (a) \\
    \includegraphics[width=0.45\textwidth]{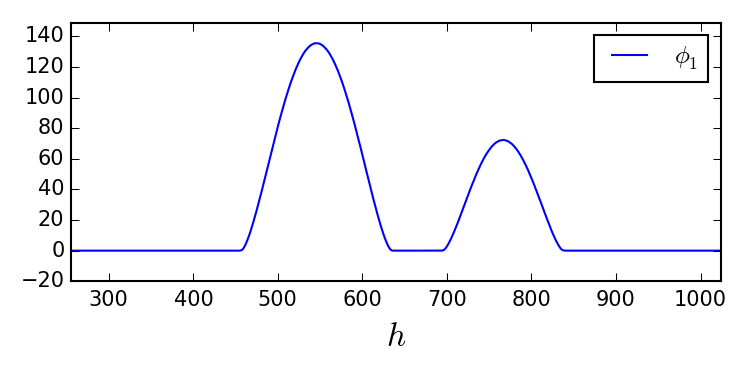} 
    &
    \includegraphics[width=0.40\textwidth]{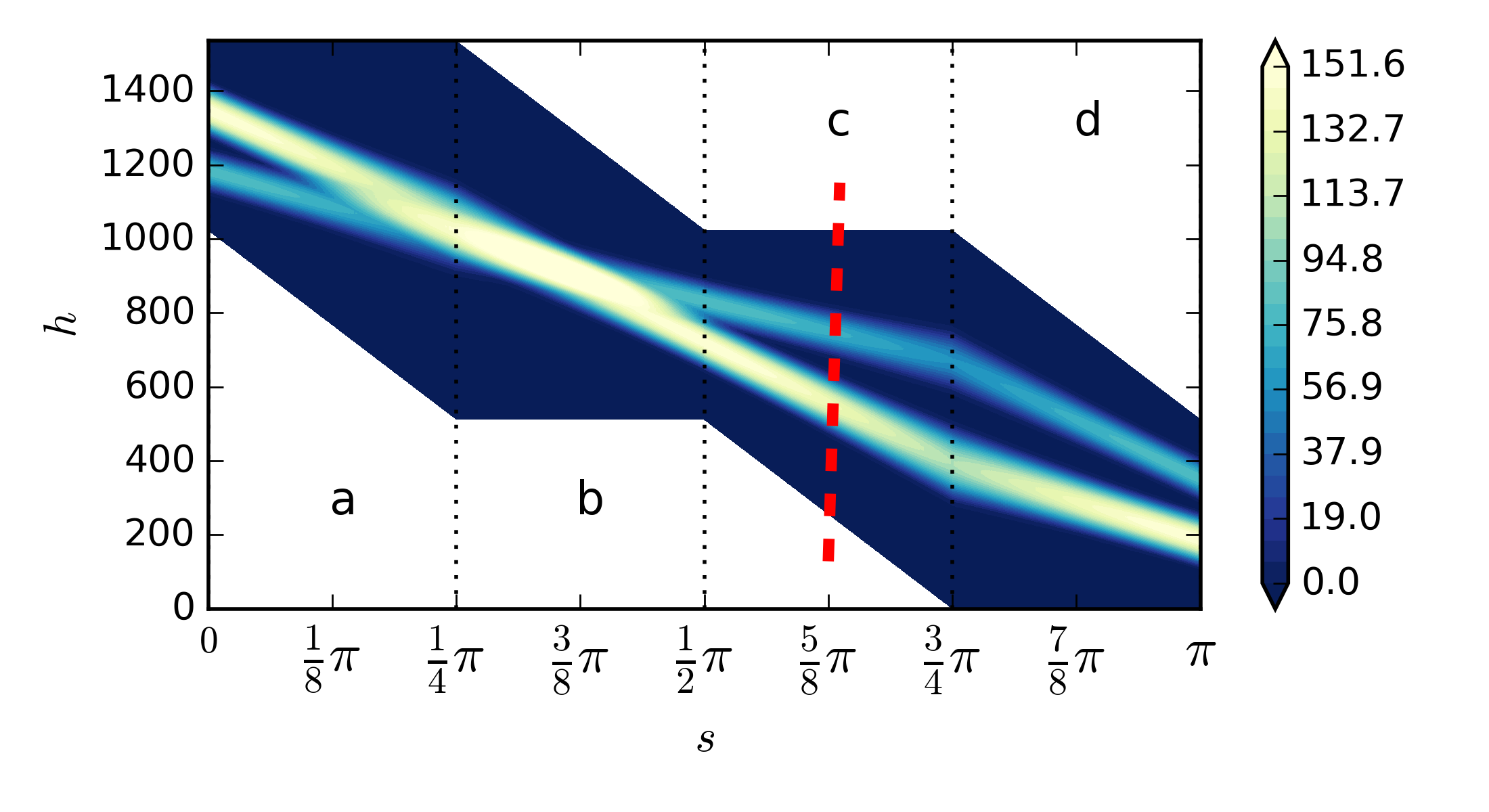} 
    \\
        (b) \\
    \includegraphics[width=0.45\textwidth]{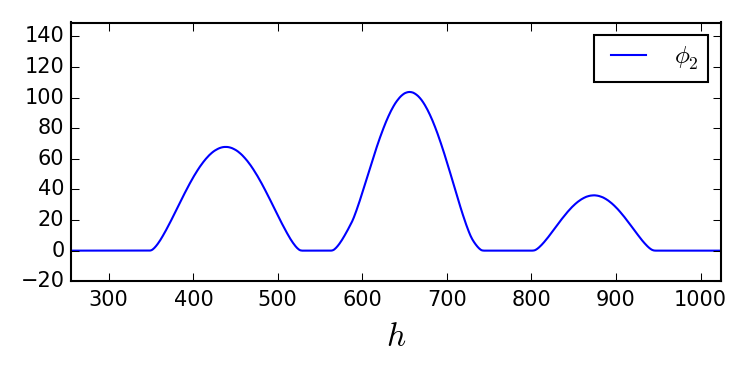} 
    &
    \includegraphics[width=0.40\textwidth]{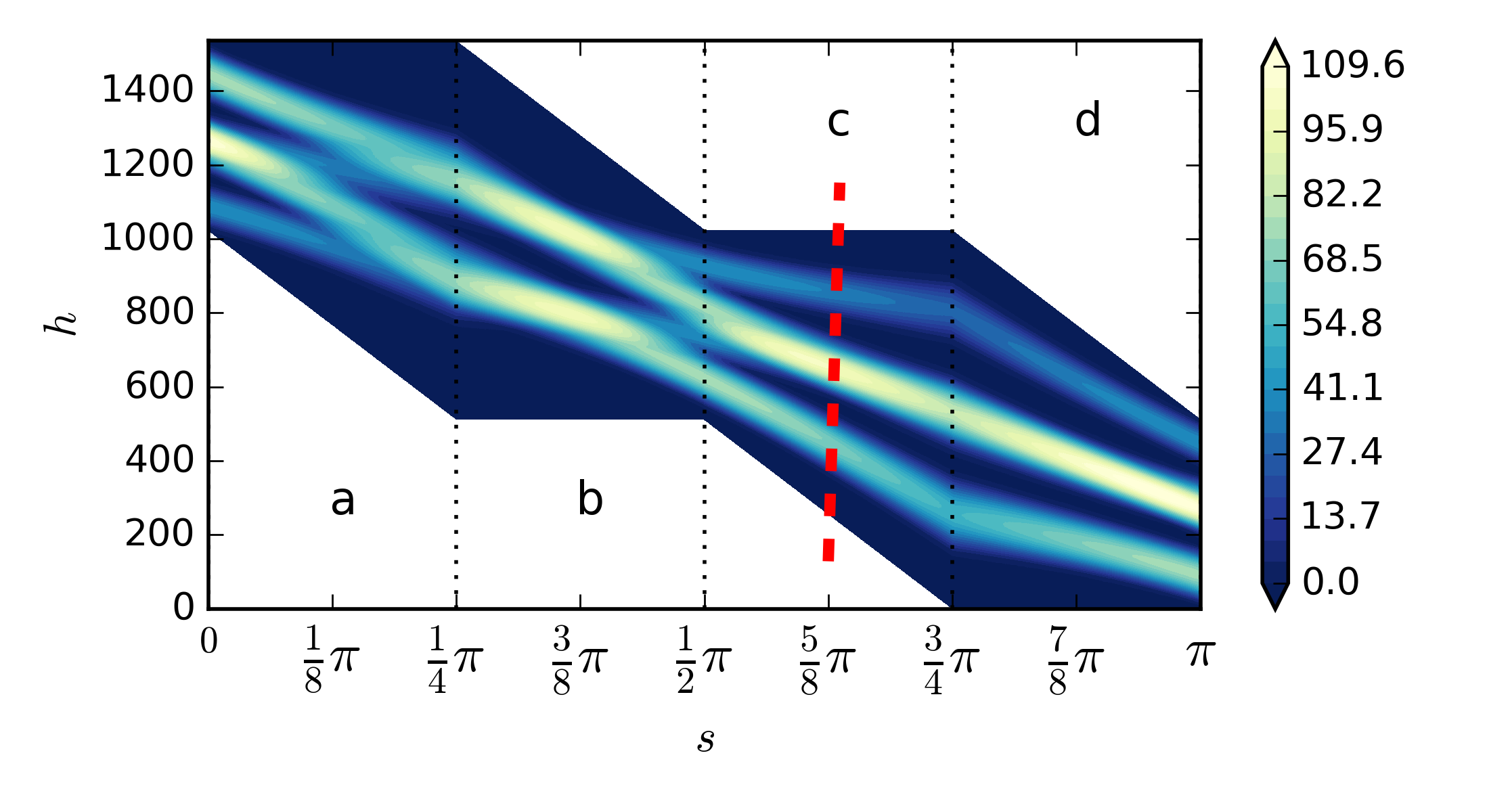} 
    \\
        (c) \\
    \includegraphics[width=0.45\textwidth]{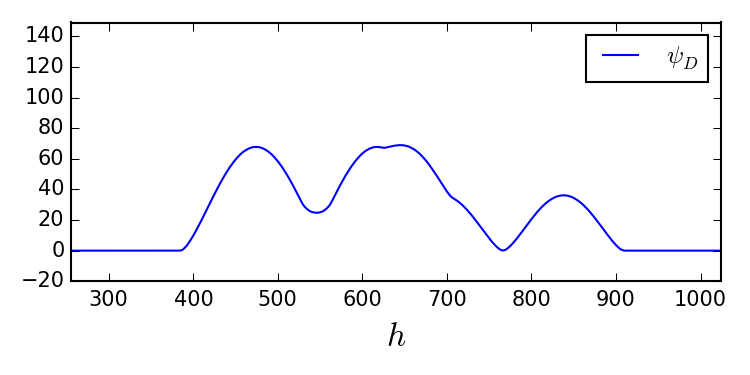} 
    &
    \includegraphics[width=0.40\textwidth]{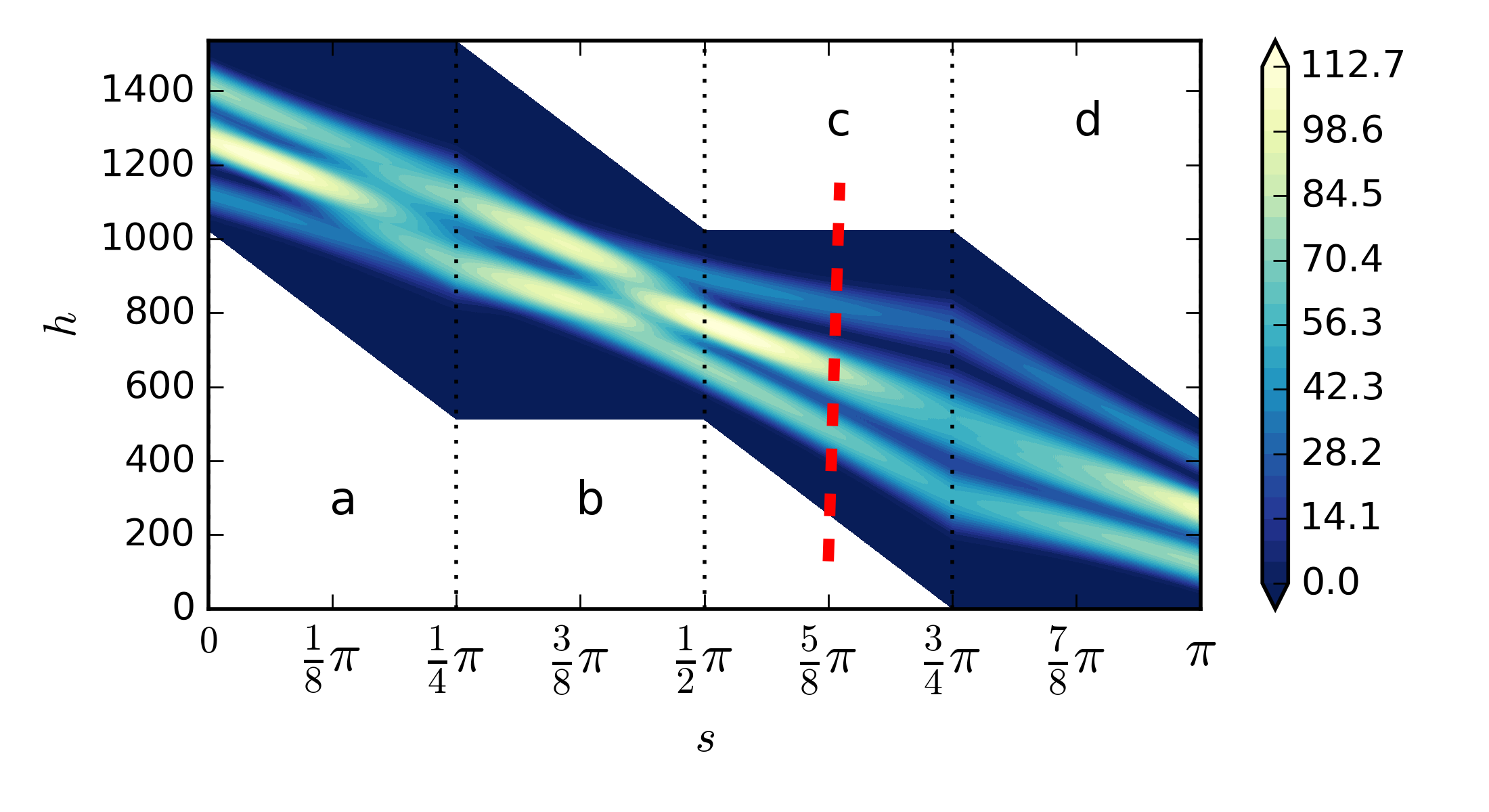} 
    \end{tabular}
    \caption{Two 1D functions $\phi_1$, $\phi_2$ and the displacement interpolation
        $\psi_D$  are shown in the first column. These are exactly the
        $s=\tan(\frac{5}{8}\pi)(N-1)$ slice
    of the DRT of the acoustic equation example  in Figure 
    \ref{fig:acoustic_2humps_drt}. The slice is indicated by the dashed red
    vertical line in the plots in the right column.}
    \label{fig:DRTdisp}
\end{figure}

\begin{figure}
    \begin{tabular}{cc}
    \includegraphics[width=0.45\textwidth]{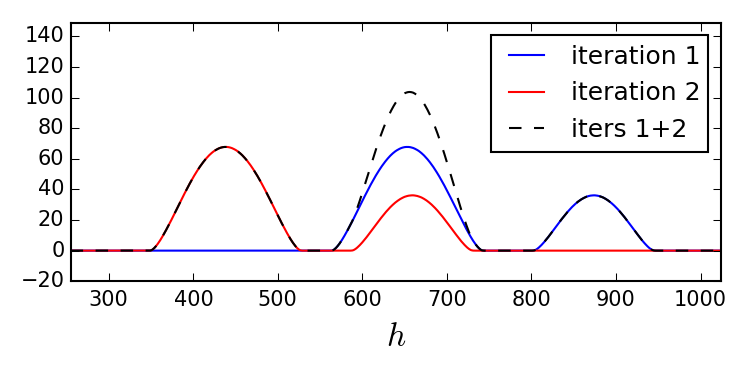} &
    \includegraphics[width=0.45\textwidth]{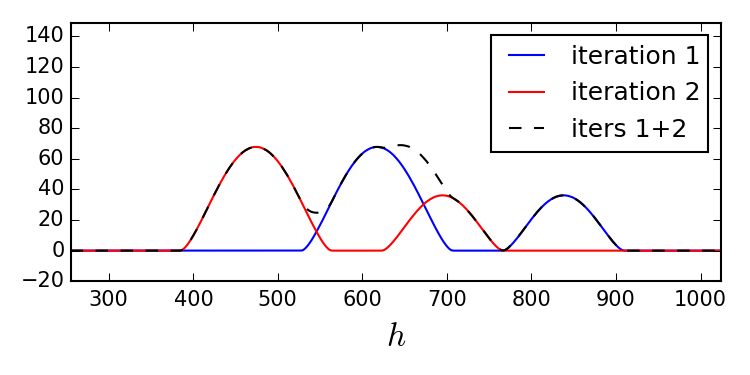} 
    \end{tabular}
    \caption{The first two contributions \eqref{eq:iterates} 
        of the transport reversal for $\phi_1$ and
    $\phi_2$ shown in Figure \ref{fig:DRTdisp} (left) and the 
    displacement interpolation resulting in $\psi_D(x;0.5)$ \eqref{eq:reversal1d}
    (right). $\psi_D$ shown in dotted line is also displayed in the bottom of 
    Figure \ref{fig:DRTdisp}.
    }
    \label{fig:tr}
\end{figure}

As an example, let us assume we are given two functions $\phi_1$ and $\phi_2$
as shown in Figure \ref{fig:DRTdisp}(a) and \ref{fig:DRTdisp}(b).  These are
taken from the 1D slice located at $s = \tan(\frac{5}{8}\pi)(N-1)$ of the DRT
from Figure \ref{fig:acoustic_2humps_drt}.  The transport reversal would
decompose $\phi_2$ into a superposition of two traveling profiles,
\begin{equation}
   \eta_1(\tau) \cK(\tau)[\rho_1(x)\varphi_1(x)]
   \quad \text{ and } \quad
   \eta_2(\tau) \cK(-\tau)[\rho_2(x)\varphi_1(x)],
   \label{eq:iterates}
\end{equation}
each plotted in Figure \ref{fig:tr}. In exact arithmetic, the two iterations of
transport reversal would pick up exactly the d'Alembert solution (although in
practice numerical error would require further iterations to pick off the
residuals). That is, we would obtain $h_1 = h_2 = 1/2$ and $\rho_1 = \rho_2 =
1$ with $\nu_1 = -\nu_2 = 1$ and $K =2$ in \eqref{eq:reversal1d}.  Now, the
displacement interpolation for $\tau = 1/2$ can be computed, yielding $\psi_D$
shown in Figure \ref{fig:DRTdisp}(c).  The exact evolution of the two iterates
\eqref{eq:iterates} are shown in Figure \ref{fig:tr}.

Now, this displacement interpolation was done for a single \emph{slice} of the
fixed $\omega$ in the transformed variables.  Suppose we are given a function
$\varphi$ in 2D.  Then by performing the same transport reversal on its Radon
transform $\hat{\varphi}$ for all $\omega$ as functions of the variable $s$,
we obtain an extension of the 1D displacement interpolant \eqref{eq:1ddi} to
higher spatial dimensions.  For each fixed angle $\omega_i \in S^1$ we obtain
the transport reversal in terms of the traveling structures,
\begin{equation}
    \hat{\psi}_D(\omega_i,s;\tau) = 
    \sum_{k=1}^K \eta_{i,k}(\tau)
    \cK(\nu_{i,k}\tau)[\rho_{i,k}(s,\tau)\hat{\varphi}(\omega_i,s)].
    \label{eq:reversal_drt}
\end{equation}
These can be used for displacement interpolation as above for each $\omega$.
The inverse transform can be taken to obtain the displacement interpolant
$\psi_D$.  

Let us clarify the implication. For the acoustic equation example with the
initial condition \eqref{eq:2humps}, we were given a snapshot of the solution
$q$ at time $t_1=0.5$ and $t_2=1.5$. From the two snapshots, we were able to
accurately approximate the solution for all time, without additional
information about the dynamics, without even knowing the PDE.  Thus this
interpolant can be more useful than linear interpolation: the linear subspace
spanned by $\{q(t_1,x),q(t_2,x)\}$ does not a contain a good approximation
for representing the evolving solution.

This ability to compute the displacement interpolation by exploiting the simple
dynamic on the transformed side will be useful in the future development of
transport reversal as a model reduction tool in multi-dimensional settings.

\section{Conclusion and future work}\label{sec:conclusion}

We have introduced a dimensional splitting method using the intertwining
property of the Radon transform. Its applications in solving hyperbolic PDEs,
imposing absorbing boundary conditions, and computing displacement
interpolations were discussed. For the inversion of DRT the conjugate gradient
method was used. 

As noted in Section \ref{sec:absorb}, the dimensional splitting proposed here
used with DRT in 3D (Section \ref{sec:drt3}) allows one to impose absorbing
boundary conditions for 3D problems without incurring any error of the type
\eqref{eq:absorberror0} that appears in 2D.  This will be verified in future
work.  The application of this splitting to fully nonlinear hyperbolic PDEs as
discussed in Section \ref{sec:nonlinearlts} will be studied as well.  The
utility of the Radon transform for displacement interpolation (Section
\ref{sec:disp}) will be much more compelling when used in conjunction with the
fully multi-dimensional transport reversal \cite{reversal}  as a model
reduction tool for general hyperbolic PDEs, and work is underway for such an
extension.

The number of CG iterations for the inversion \eqref{eq:normaleqn} can be
estimated to justify the conjectured $\cO(N^{5/2}\log N)$ cost for inversion:
this and other inversion results will appear elsewhere.  While the prologation
used in the inversion (Section \ref{sec:idrt}) causes expense only of a
constant factor, it can be of significant computational cost. Other approaches
to reduce the amount of computational effort will be explored.  The DRT is
essentially a structured matrix multiplication and may be amenable to
parallelization, and its performance on graphical processing units (GPUs) will
be a future topic of research.

\section*{Acknowledgments} 
The author expresses gratitude to Randall J. LeVeque for carefully reviewing this
manuscript. He also thanks Gunther Uhlmann for helpful discussions.

\bibliographystyle{siamplain}

\end{document}